\documentclass[review]{elsarticle}

\usepackage{lineno,hyperref}
\modulolinenumbers[5]
\journal{Elsevier}
\usepackage[a4paper, total={16cm, 21.0cm}]{geometry}
\usepackage{amsmath,graphicx, xcolor,amssymb}
\usepackage[font=small,labelfont=bf,tableposition=top]{caption}
\usepackage{graphicx}
\usepackage{epstopdf}
\usepackage{mathabx}
\usepackage{mathrsfs}
\usepackage{ragged2e}
\usepackage[mathscr]{euscript}
\usepackage{multirow}
\usepackage{subfigure}
\usepackage{enumerate}
\usepackage{mathrsfs}
\usepackage{color,soul}
\usepackage{todonotes}
\usepackage{verbatim}
\usepackage{float}
\usepackage{framed}
\usepackage{pagecolor}
\usepackage{algorithm, algorithmic}
\graphicspath{{figures/}{figs/}{figures/Exp1/}{figures/Exp2/}{figures/Exp3/}{fig/}}

\begin{document}

\begin{frontmatter}


\title{A multiscale Robin-coupled implicit method for two-phase flows in high-contrast formations}


\author[mymainaddress]{Franciane F. Rocha}
\ead{fr.franciane@usp.br}
\corref{mycorrespondingauthor}
\cortext[mycorrespondingauthor]{Corresponding author}

\author[mymainaddress]{Fabricio S. Sousa}
\ead{fsimeoni@icmc.usp.br}

\author[mymainaddress]{Roberto F. Ausas}
\ead{rfausas@icmc.usp.br}

\author[mymainaddress]{Gustavo C. Buscaglia}
\ead{gustavo.buscaglia@icmc.usp.br}

\author[mysecondaryaddress]{Felipe Pereira}
\ead{luisfelipe.pereira@utdallas.edu}

\address[mymainaddress]{Instituto de Ci\^encias Matem\'aticas e de Computa\c c\~ao, Universidade de S\~ao Paulo,\\ Av. Trabalhador S\~ao-carlense, 400, 13566-590, S\~ao Carlos, SP, Brazil}
\address[mysecondaryaddress]{Department of Mathematical Sciences, The University of Texas at Dallas,\\ 800 W. Campbell Road, Richardson, TX 75080-3021, USA}

\begin{abstract}
In the presence of strong heterogeneities, it is well known that the use of explicit schemes for the transport of species in a porous medium suffers from severe restrictions on the time step. This has led to the development of implicit schemes that are increasingly favoured by practitioners for their computational efficiency. The transport equation requires knowledge of the velocity field, which results from an elliptic problem (Darcy problem) that is the most expensive part of the computation. When considering large reservoirs, a cost-effective way of approximating the Darcy problems is using multiscale domain decomposition (MDD) methods. They allow for the pressure and velocity fields to be computed on coarse meshes (large scale), while detailed basis functions are defined locally, usually in parallel, in a much finer grid (small scale). In this work we adopt the Multiscale Robin Coupled Method (MRCM, [Guiraldello, et al., J. Comput. Phys., 355 (2018) pp. 1-21], [Rocha, et al., J. Comput. Phys., (2020) 109316]), which is a generalization of previous MDD methods that allows for great flexibility in the choice of interface spaces. In this article we investigate the combination of the MRCM with implicit transport schemes. A sequentially implicit strategy is proposed, with different trust-region algorithms ensuring the convergence of the transport solver. The method is assessed on several very stringent 2D two-phase problems, demonstrating its stability even for large time steps. It is also shown that the best accuracy is achieved by considering recently introduced non-polynomial interface spaces, since polynomial spaces are not optimal for high-contrast channelized permeability fields. 

\end{abstract}

\begin{keyword}
Multiscale Robin Coupled Method \sep coupled flow and transport \sep sequential implicit solution \sep high-contrast porous media
\end{keyword}

\end{frontmatter}



\section{Introduction}

Multiscale domain decomposition methods are a suitable choice to deal with the huge elliptic problems arising from the discretization of the equations governing multiphase flows in oil reservoirs \cite{chen2006computational}. They allow for the pressure and velocity fields to be computed on a coarse mesh (large scale), while detailed basis functions (locally defined for each subdomain) incorporate rock heterogeneity on a much finer grid (small scale) \cite{kippe}. The local problems can be solved simultaneously in state-of-the-art parallel machines, making the simulation of huge problems feasible \cite{abreu2020recursive}. 

Several multiscale methods have been presented in the context of the finite volume method \cite{jenny2003multi,jenny2005adaptive}, the finite element method \cite{HouMulti, aarnes2002multiscale, efendiev2013generalized}, and mixed finite elements \cite{chen2003mixed, aarnes, arbogast, pereira, chung2015mixed}. We consider here the Multiscale Robin Coupled Method (MRCM \cite{guiraldello2018multiscale, bifasico}) for the solution of two-phase flow problems. It is a generalization of the Multiscale Mixed Method (MuMM \cite{pereira}) that allows for the independent choice of pressure and flux interface spaces through local Robin boundary conditions. For improved accuracy when the permeability field has high contrast and is channelized, we incorporate adapted interface spaces which  behave better than classical polynomials \cite{guiraldello2019interface, rocha2020interface}.

Our focus in this work is the solution of nonlinear two-phase flow models. 
In the literature the coupling of multiscale flow and transport problems has been treated by explicit operator splitting techniques  \cite{furtado2011operator, paz2020adaptive} and implicit formulations \cite{jenny2006adaptive, ganis2014global}. We have considered operator splitting techniques with explicit approximations for the transport problem in previous works \cite{bifasico, rocha2020interface}. Here, we propose to combine the MRCM with a Sequential Implicit (SI) scheme \cite{watts1986compositional} that allows for the use of large time steps for the coupled flow and transport problem, in contrast to explicit time integration approaches where a CFL-type condition restricts the size of the time step for the transport calculation. 
The SI algorithm solves at each time step one update of the pressure equation followed by the implicit solution of the transport by the Newton method, considering the velocity field fixed in time. 

To ensure convergence of the nonlinear loop we consider trust-region methods to guide the Newton iterations. Specifically, our new solver has three options of trust-region algorithms. The first is the strategy of Jenny et al. \cite{jenny2009unconditionally}, where two successive iterations cannot cross any trust-region boundary delineated by the inflection-point of the analytic flux function. The other two algorithms are based on least-square methods, that automatically select the trust-regions to define the iterative step of the implicit solver which are the trust-region reflective and the trust-region dogleg algorithms \cite{byrd1988approximate, powell1968fortran, conn2000trust}. 
We remark that many improvements have been proposed to the inflection-point strategy \cite{jenny2009unconditionally}, for example, the extension to flows with buoyancy and capillary forces \cite{wang2013trust}, compositional simulations \cite{voskov2011compositional}, and the development of a numerical trust-region solver based on the discretized flux function \cite{li2015nonlinear}. 
Other developments based on trust-region methods are available in the literature. For instance, in \cite{moyner2017nonlinear} the authors developed a flux-search solver for three-phase flow problems. 

In this paper we investigate the behavior of domain-decomposition-based multiscale methods (specifically, the family of methods parameterized by the MRCM) when coupled with implicit transport solvers for the simulation of two-phase flows. The approximation is assessed in terms of accuracy and computational efficiency. Additionally, we compare three nonlinear iterative schemes considering different trust-region algorithms so as to determine the one that provides the best performance for the proposed multiscale SI method.

Though satisfactory results in the studied flows were obtained, the semi-implicit treatment of the velocity can generate material balance errors in more complex models (e.g., compositional flows)  \cite{aziz1979petroleum}. The Sequential Fully-Implicit (SFI) scheme, developed in the Multiscale Finite Volume Method (MSFV) framework \cite{jenny2006adaptive}, is an option to deal with this issue. See  \cite{moncorge2017modified, moncorge2018sequential,  moyner2020nonlinear,  lee2020conservative} for extensions to compositional flow simulations, and \cite{jiang2019nonlinear} for developments on nonlinear acceleration techniques. 
We also study in this work the MRCM combined with the SFI scheme to approximate a challenging two-phase flow problem, with strong fingering instabilities \cite{glimm1981numerical}. 

Summarizing, the main contributions of this work are:
\begin{itemize}
\item[•] We show that the family of multiscale methods considered here can be efficiently coupled, through both SI and SFI schemes, with implicit transport solvers for the simulation of two-phase flows. 
\item[•] We identify the trust-region algorithm, with the inflection-point strategy, as the best-performing nonlinear iterative method for the aforementioned coupling.
\end{itemize}



The rest of the paper is organized as follows: In Section \ref{sec_Two}, the two-phase flow model is presented. We recall the MRCM in Section \ref{sec_MRCM} and present the sequential implicit formulation in Section \ref{sec_SFI}. 
Numerical simulation results are presented in Section \ref{sec_results}. Finally, our conclusions are presented in Section \ref{sec_conclusions}.

\section{Two-phase flows}\label{sec_Two}

The governing system of equations for the two-phase problem consists of a second-order elliptic equation for pressure coupled to a hyperbolic conservation law for the saturation of one of the phases \cite{ewing1983mathematics}.
The phases considered are water and oil (denoted by $w$ and $o$, respectively), and the sum of their saturations is equal to one since we  assume a fully saturated porous medium. 
We consider an immiscible and incompressible two-phase flow in a reservoir containing injection and production wells. For simplicity, capillary pressure effects are not taken into account. However, the proposed method could incorporate these effects by means of an operator splitting algorithm \cite{douglas1997numerical}.
Adopting a mixed formulation for the pressure equation and a finite volume scheme for the hyperbolic conservation law, the unknowns of the two-phase flow problem are the Darcy velocity $\mathbf{u}(\mathbf x, t)$, the fluid pressure $p(\mathbf x, t)$ and the water saturation $s(\mathbf{x},t)$. The pressure and velocity are related by Darcy's law so that the elliptic problem can be written as
\begin{equation}\label{Darcy}
\begin{array}{rll}
\mathbf{u}&=-\lambda(s)K(\mathbf{x})\nabla p  &\mbox{in}\ \Omega \\
\nabla \cdot \mathbf{u}&=q  &\mbox{in}\ \Omega \\ 
p &= p_b &\mbox{on}\ \partial\Omega_{p}\\
\mathbf{u} \cdot \mathbf{n}&= u_b, &\mbox{on}\ \partial\Omega_{u}
\end{array}
\end{equation}
where $\Omega\subset\mathbb{R}^d,\ d=2$ or $d=3$ is the domain, $\partial\Omega=\partial\Omega_{p}\cup\partial\Omega_{u}$, $\partial\Omega_{p}\cap\partial\Omega_{u}=\emptyset$; $K(\mathbf{x})$ is the symmetric, uniformly positive definite absolute permeability tensor; $q = q(\mathbf x, t)$ is the source term; $p_b = p_b(\mathbf x, t)$ is the pressure boundary condition at boundary $\partial\Omega_{p}$; $u_b = u_b(\mathbf x, t)$ is the normal velocity boundary condition ($\mathbf{n}$ is the outward unit normal) at the boundary $\partial\Omega_{u}$; $\lambda(s)$ is the total mobility, given by the sum of the mobilities of the phases:
\begin{equation}
 \lambda(s) =\lambda_w(s)+\lambda_o(s)=\dfrac{k_{rw}(s)}{\mu_w}+\dfrac{k_{ro}(s)}{\mu_o},
 \label{fluxos1} 
\end{equation}
where $k_{rj}(s)$ and $\mu_j$, $j \in \{w,o\} $, are, respectively, the relative permeability and viscosity of phase $j$. 
The water saturation problem is governed by the transport equation
\begin{equation}\label{BL2D}
\begin{array}{rll}
\dfrac{\partial s}{\partial t} + \nabla \cdot \left(f(s)\mathbf{u}\right)& = 0  &\mbox{in}\ \Omega \\
s(\mathbf{x},t=0) &= s^0(\mathbf{x}) &\mbox{in}\ \Omega\\
s(\mathbf{x},t) &= \bar{s}(\mathbf{x},t) &\mbox{in}\ \partial\Omega^-
\end{array}
\end{equation}
where $s^0$ and $\bar{s}$ are, respectively, the initial and injection conditions for the water saturation. Here, $\partial\Omega^-=\{\mathbf{x}\in\partial\Omega,\ \mathbf{u} \cdot \mathbf{n}<0\}$ represents the inlet boundaries. The function $f(s)$ is the nonlinear fractional flow of water, given by
\begin{equation}
\qquad f(s) = \dfrac{\lambda_w(s)}{\lambda(s)}.
 \label{fluxos2} 
\end{equation}
For simplicity, we assume a constant porosity (scaled out by changing the time variable). 

We combine the MRCM to solve equation (\ref{Darcy}) with an implicit approximation of the hyperbolic conservation law for the water saturation (\ref{BL2D}) in a sequential fashion, aiming at a compromise between accuracy and computational efficiency of numerical simulations.

\section{The Multiscale Robin Coupled Method}\label{sec_MRCM}

The Multiscale Robin Coupled Method has been introduced to solve elliptic equations accurately, presenting advantages when compared to other existing multiscale mixed methods for flows in highly heterogeneous porous media \cite{guiraldello2018multiscale, bifasico}. 
The MRCM solves equation (\ref{Darcy}) by a decomposition of the domain $\Omega$ into non-overlapping subdomains $\Omega_i, \ i=1,2,\cdots,N$. Weak continuity of normal fluxes and pressures are imposed to the multiscale solution $(\mathbf u_h,p_h)$ at the skeleton $\Gamma$ of the decomposition (the union of the interfaces $\Gamma_{ij} = \Omega_i \cap \Omega_j$), whose characteristic size $H$  is significantly larger than the fine-scale of the discretization ($H\gg h$).
The weak continuities are enforced by the following compatibility conditions 
\begin{equation}
\int_\Gamma (\mathbf u_h^+ - \mathbf u_h^-)\cdot \check{\mathbf n} \ \psi \ d\Gamma=0 \quad \text{and}\quad 
\int_\Gamma (p_h^+ - p_h^-)\ \phi \ d\Gamma=0
\label{eq:compatibility}
\end{equation}
for all $(\phi,\psi) \in \mathcal{U}_H \times \mathcal{P}_H$, that are the interface spaces defined over the edges $\mathcal{E}_h$ of the skeleton $\Gamma$ as subspaces of
\begin{equation}
\mathfrak{F}_h(\mathcal{E}_h) = \left\{ f:\mathcal{E}_h\to \mathbb{R} ; ~f|_e\,\in\,\mathbb{P}_0~,~\forall\,e\,\in\,\mathcal{E}_h \right\} ~.
\end{equation}
In Eq. (\ref{eq:compatibility}), the $+$ and $-$ superscripts represent the solution on each side of the interface $\Gamma$, while $\check{\mathbf n}$ is a fixed normal vector to the skeleton $\Gamma$ (pointing outwards from the subdomain with the smallest index).

The formulation of the MRCM consists of finding local solutions $(\mathbf u_h^i, p_h^i)$ within each subdomain $\Omega_i$, and interface unknowns $(U_H,P_H)$ satisfying 
\begin{equation}
\begin{array}{rcll}
\mathbf u_h^i &=& -\kappa \ \nabla p_h^i & \text{in }\Omega_i \\
\nabla \cdot \mathbf u_h^i &=& q & \text{in }\Omega_i \\
p_h^i & = & p_b & \text{on }\partial\Omega_i\cap \partial\Omega_p \\
\mathbf u_h^i \cdot \check{\mathbf n}^i &=& u_b & \text{on }\partial\Omega_i\cap \partial\Omega_u \\
-\beta_i  \mathbf u_h^i \cdot \check{\mathbf n}^i + p_h^i &=&  -\beta_i U_H \check{\mathbf n} \cdot \check{\mathbf n}^i + P_H & \text{on }\partial\Omega_i\cap \Gamma 
\end{array}
\label{eq:mrc1}
\end{equation}
along with the compatibility conditions
\begin{equation}
\begin{array}{rcl}
\displaystyle \sum_{i=1}^N \int_{\partial\Omega_i\cap \Gamma} (\mathbf u_h^i \cdot \check{\mathbf n}^i) \ \psi \ d\Gamma &=& 0\\
\displaystyle \sum_{i=1}^N \int_{\partial\Omega_i\cap \Gamma} \beta_i (\mathbf u_h^i \cdot \check{\mathbf n}^i - U_H\ \check{\mathbf n} \cdot \check{\mathbf n}^i)\ \phi  \ (\check{\mathbf n} \cdot \check{\mathbf n}^i)\ d\Gamma &=& 0
\end{array}
\label{eq:mrc2}
\end{equation}
for all $(\phi,\psi)\in \mathcal{U}_H\times\mathcal{P}_H$, where $\kappa = \lambda(s(\mathbf x)) K(\mathbf x)$, and $\check{\mathbf n}^i$ is the normal vector to $\Gamma$ pointing outwards of $\Omega_i$.
The parameter for the local Robin boundary conditions can be written as
\begin{equation}
\beta_i(\mathbf {x}) = \frac{\alpha(\mathbf{x}) H}{\kappa_i(\mathbf x)},
\label{eq:beta}
\end{equation}
where $\alpha(\mathbf {x})$ is a dimensionless algorithmic function. We remark that the variation of this dimensionless parameter, along with a suitable choice of interface spaces, result in a family of different methods, and by setting this function to extreme values ($\alpha \rightarrow 0$ and $\alpha \rightarrow +\infty$) well-known multiscale mixed methods can be recovered \cite{guiraldello2018multiscale}. 

The implementation of the MRCM consists of solving the equations (\ref{eq:mrc1}) independently to obtain a local set of multiscale basis functions for each subdomain. The global solution is then given by a linear combination of the basis functions and coefficients obtained from the global interface system (\ref{eq:mrc2}).

\subsection{Interface spaces}
It is well known that the classical polynomials are not optimal for high-contrast channelized permeability fields \cite{guiraldello2019interface}. Here we recall novel interface spaces based on physics to deal with permeability fields containing highly-permeable channels and barriers, that has been introduced recently \cite{rocha2020interface, rochaenhanced}. These interface spaces are particularly relevant when the channels and barriers are relatively large as happens in karst reservoirs \cite{ popov2009multiphysics,  lopes2019new}.

Let $\Gamma_{i,j}\subset\Gamma$ be an interface with support in the line segment $[a,b]$ through which $N_{\text{high}}$ high-permeability channels pass. Let $[a_k,b_k]\subset [a,b]$, $k=1,\cdots,N_{\text{high}}$ denote the respective support of each channel. 
The pressure space contains the following basis functions that mimic the behavior of the pressure solution across the channels:  
\begin{equation}\label{new_basis1}
\psi_0(x)=\left\{
\begin{array}{cl}
\dfrac{a_1-x}{a_1-a}  &\mbox{if}\ x\in(a,a_1) \\
0\   &\mbox{otherwise}
\end{array} \right. 
\end{equation}
\begin{equation}\label{new_basis3}
\psi_{N_{\text{high}}+1}(x)=\left\{
\begin{array}{cl}
\dfrac{x-b_{N_{\text{high}}}}{b-b_{N_{\text{high}}}}  &\mbox{if}\ x\in(b_{N_{\text{high}}},b) \\
0\   &\mbox{otherwise}
\end{array} \right.
\end{equation}

\begin{equation}\label{new_basis2}
\psi_k(x)=\left\{
\begin{array}{cl}
\dfrac{x-b_{k-1}}{a_{k}-b_{k-1}}  &\mbox{if}\ x\in(b_{k-1},a_{k}) \\
1  &\mbox{if}\ x\in(a_{k},b_{k}) \\
\dfrac{a_{k+1}-x}{a_{k+1}-b_{k}}  &\mbox{if}\ x\in(b_{k},a_{k+1}) 
\end{array} \right.
\end{equation}
for $k=1,\cdots,N_{\text{high}}$. We assume $b_0=a$ and $a_{N_{\text{high}}+1}=b$ in Eq. (\ref{new_basis2}). Therefore, the total number of pressure basis functions at $\Gamma_{i,j}$ is $2+N_{\text{high}}$.

Next, let $\Gamma_{i,j}\subset\Gamma$ be an interface with support in $[a,b]$ through which $N_{\text{low}}$ low-permeability structures pass. Let $[a_k,b_k]\subset [a,b]$, $k=1,\cdots,N_{\text{low}}$ the respective support of each low-permeability structure. The flux space mimics the behavior of the flux across the barriers and contains the following basis functions:

\begin{equation}\label{new_basis4}
\phi_0(x)=\left\{
\begin{array}{cl}
1  &\mbox{if}\ x\in(a,a_1) \\
0\   &\mbox{otherwise}
\end{array} \right.
\end{equation}
\begin{equation}\label{new_basis6}
\phi_{2k-1}(x)=\left\{
\begin{array}{cl}
1  &\mbox{if}\ x\in(a_k,b_k) \\
0\   &\mbox{otherwise}
\end{array} \right.
\end{equation}

\begin{equation}\label{new_basis5}
\phi_{2k}(x)=\left\{
\begin{array}{cl}
1  &\mbox{if}\ x\in(b_k,a_{k+1}) \\
0\   &\mbox{otherwise} 
\end{array} \right.
\end{equation}
for each low-permeability structure $k=1,\cdots,N_{\text{low}}$, where we assume $a_{N_{\text{low}}+1}=b$. Therefore the total of flux basis functions at $\Gamma_{i,j}$ is $1+2N_{\text{low}}$.

We use the MRCM with the physics-based interface spaces to capture the geometry of the high-permeability channels and low-permeability structures at each interface. We also consider the adaptivity in the $\alpha(\mathbf {x})$ function according to permeability variations \cite{bifasico}. Thus, one can control the relative importance of each interface space at each location. As proposed in \cite{bifasico}, we take a small value (pressure is favored) at the high-permeability channels and a large value (flux is favored) for the remaining areas.
The interface spaces at the interfaces without high-permeability channels or low-permeability structures are chosen as linear polynomials.

\section{A sequential implicit solver for two-phase subsurface flows} \label{sec_SFI}

In the presence of strong heterogeneity, explicit schemes for the transport of saturation suffer from severe time-step restrictions. To solve the coupled equations (\ref{Darcy}) and (\ref{BL2D}) we consider the sequential implicit method \cite{watts1986compositional}, which allows for the use of large time steps, improving the computational efficiency of the simulation.

In the SI algorithm, each time step consists of a sequential update for the flow and transport problems, where a (nonlinear) Newton loop is used to solve the transport equation implicitly. 
We denote by $\Delta t$ the time step used to update the coupled problems of flow and transport at times $t^n = n\Delta t$, for $n=0,1,\dots$. 
Let $p^n (\mathbf x)$, $\mathbf u^n (\mathbf x)$ and $s^n (\mathbf x)$ denote the pressure, velocity and saturation approximations for $p(\mathbf x, t^n)$,  $\mathbf u(\mathbf x, t^n)$ and $s(\mathbf x, t^n)$ respectively, at time $t^n$.
To find the updated variables, one first computes $s^{n+1} (\mathbf x)$ (as detailed below) and then solves  (\ref{Darcy}) for the pressure $p^{n+1} (\mathbf x)$ and velocity $\mathbf u^{n+1} (\mathbf x)$ keeping the saturation frozen at $s^{n+1}$.

The saturation $s^{n+1} (\mathbf x)$ is computed through Eq. (\ref{BL2D}) by using a simple implicit Euler time integration considering $\mathbf{u}$ constant in time as follows 
\begin{equation}\label{BL2D_Newton}
\dfrac{s^{n+1}-s^{n}}{\Delta t} + \nabla \cdot \left(f(s^{n+1})\mathbf{u}^n\right) = 0.  
\end{equation}
Upon a finite volume discretization, the problem for $s^{n+1}$ reads
\begin{equation}\label{BL2D_Newton2}
s^{n+1}_{I}=s^{n}_{I} - \dfrac{\Delta t}{V_I} \Big(\mathcal{F}_I^{n+1} \Big),
\end{equation}
where $I$ refers to a computational cell of an orthogonal, uniformly spaced (by directions) grid identified as an index ($I=(i,j)$ in 2D and $I=(i,j,k)$ in 3D), $ V_I$ represents the volume of cell $I$, and $\mathcal{F}_I^{n+1}$ is a function of $f(s^{n+1,\nu+1})$ and $\mathbf{u}^n$, that represents the balance of fluxes at the faces of cell $I$. 

We solve Eq. (\ref{BL2D_Newton2}) by variants of Newton's iterative method. Let $\nu$ refer to the iteration level of the Newton loop for saturation and set $s^{n+1,0}=s^n$.
In a pure Newton scheme, the next iterate $s^{n+1,\nu+1}$ is defined by the linear system
\begin{equation}\label{BL2D_Newton8}
    \mathcal{H}'(s^{n+1,\nu})\ {\bf d}^{\nu} =-\mathcal{H}(s^{n+1,\nu}),
\end{equation}
where
\begin{equation}\label{BL2D_Newton7}
    \mathcal{H}(s^{n+1,\nu})= \left[s^{n+1,\nu} - s^{n} + \dfrac{\Delta t}{V}\mathcal{F}^{n+1,\nu}\right]_I,
\end{equation}
$\mathcal{H}'(s^{n+1,\nu})$ is the Jacobian matrix of $\mathcal{H}$ and $s^{n+1,\nu+1}=s^{n+1,\nu}+{\bf d}^{\nu}$. 
The solution at the new time level is achieved when the change in the saturation between two successive iterations is less than a specified tolerance denoted by $\eta$. In other words,  $s^{n+1}=s^{n+1,\nu+1}$ if $\parallel s^{n+1,\nu+1}-s^{n+1,\nu}\parallel\leq \eta$.

We consider here a first-order upwind scheme to define $\mathcal{F}^{n+1,\nu}_{I}$. In the 2D case we have 
\begin{equation}
\mathcal{F}^{n+1,\nu}_{I}=\mathcal{F}^{n+1,\nu}_{i,j}=\Big(F^{n+1,\nu}_{i+1/2,j} - F^{n+1,\nu}_{i-1/2,j}\Big)+  \Big(G^{n+1,\nu}_{i,j+1/2} - G^{n+1,\nu}_{i,j-1/2}\Big),
\label{fluxes_impli}
\end{equation}
with discrete fluxes $F^{n+1,\nu}_{i-1/2,j}$ and $G^{n+1,\nu}_{i,j-1/2}$ on respective interfaces $x_{i-1/2}$ and $y_{j-1/2}$ given by
\begin{equation}
 F^{n+1,\nu}_{i-1/2,j} =
 \left\{
\begin{array}{rl}
\Delta y\ f^{n+1,\nu}_{i-1,j}u^x_{i-1/2,j} &\mbox{if}\ u^x_{i-1/2,j}>0 \\
\Delta y\ f^{n+1,\nu}_{i,j}u^x_{i-1/2,j} &\mbox{otherwise}
\end{array} \right.
 \label{up_fluxes_impli}
\end{equation}
and
\begin{equation}
 G^{n+1,\nu}_{i,j-1/2} =
 \left\{
\begin{array}{rl}
\Delta x\ f^{n+1,\nu}_{i,j-1}u^y_{i,j-1/2} &\mbox{if}\ u^y_{i,j-1/2}>0 \\
\Delta x\ f^{n+1,\nu}_{i,j}u^y_{i,j-1/2} &\mbox{otherwise}
\end{array} \right. 
 \label{eq:up_fluxes_impli}
\end{equation}
where $ u^x=\  u^x(x,y)$ and $u^y= \ u^ y(x,y)$ denote the $x$ and $y$ components of the velocity field $\textbf{u}$, and $f_{i,j}^{n+1,\nu}=f(s_{i,j}^{n+1,\nu})$. The variable $s^{n+1,\nu}_{i,j}=s(x_i,y_j,t^{n+1})$ represents the saturation (assumed to be a piecewise constant over each computational cell) at time $t=t^{n+1}$ and at Newton iteration $\nu$. 

To ensure the convergence of the nonlinear loop we consider trust-region algorithms instead of pure Newton iterations. In the next subsection we present more details about the algorithms that we have implemented and compared.

\subsection{Trust-region algorithms}
\label{subsec:Trust-region algorithms}

Nonlinearity is a challenging issue for reservoir simulations. 
Convergence failures for the transport problem are related to the nonlinearity of the flux function. The Newton method is not guaranteed to converge for large time steps, and it can be sensitive to the initial guess \cite{jenny2009unconditionally}. 
To ensure convergence of the nonlinear loop we consider trust-region algorithms to guide the Newton iterations. 
Specifically, our solver has three options of trust-region algorithms: the {\em inflection-point strategy} of Jenny et al. \cite{jenny2009unconditionally}, the trust-region {\em reflective algorithm} \cite{byrd1988approximate}, and the trust-region {\em dogleg algorithm} \cite{powell1968fortran}. 

The inflection-point strategy was introduced to deal with Newton's initial guesses that are on the opposite side of the saturation inflection point with respect to the saturation solution \cite{jenny2009unconditionally}. To ensure the convergence of the Newton iterative process for any time step size, two successive saturation updates are made on the same side of the saturation inflection-point. Hence, if an update would cross the inflection point, selective under-relaxation is applied, i.e., if $f''(s^{n+1,\nu+1})f''(s^{n+1,\nu})<0$, then $s^{n+1,\nu+1}=(s^{n+1,\nu+1}+s^{n+1,\nu})/2$. Additionally, it is necessary to enforce the constraint $0\leq s^{n+1,\nu+1} \leq 1$
after every iteration, that is justified by the physics of the problem \cite{jenny2009unconditionally}. 

The inflection-point strategy can be seen as a trust-region method that defines different saturation regions delineated by the inflection-point. The updates are performed such that two successive iterations cannot cross any trust-region boundary. A more general trust-region Newton method, that includes saturation trust-regions delineated by the unit-flux and endpoints, is presented in \cite{wang2013trust}. Another extension of the inflection-point strategy was the development of a numerical trust-region solver, that is based on the discretized flux function \cite{li2015nonlinear}. However, for the problem at hand (two-phase flows without gravity and capillary effects) the inflection-point strategy is a particular case of the methods of \cite{wang2013trust} and \cite{li2015nonlinear}.




The other two trust-region algorithms considered here are quite popular for nonlinear least-squares problems \cite{nocedal2006numerical}. They define the iterative update $\mathbf{d}^{\nu}$ by minimizing a model function in a selected region \cite{yuan2015recent}. 
To explain this approach, consider the unconstrained minimization problem
\begin{equation}
\min_{\mathbf{x}\in\mathbb{R}^n}\varphi(\mathbf{x}),
\end{equation}
where $\mathbf{x}$ represents the vector with unknowns $s_{i,j}^{n+1}$ and $\varphi:\mathbb{R}^n\to\mathbb{R}$ is the objective function to be minimized, i.e.,
\begin{equation}
    \varphi(\mathbf{x})=\|\mathcal{H}(\mathbf{x})\|_2^2~.
\end{equation}
The trust-region algorithm at iteration $\nu$ defines $\mathbf{d}^{\nu}$ by solving the following sub-problem
\begin{equation}
\min_{\mathbf{d}}\{m^\nu(\mathbf{d});\ \parallel \mathbf{d} \parallel\leq \Delta^\nu,\ \mathbf{d}\in\mathbb{R}^n\}, 
\end{equation}
where $m^\nu$ is a model function that represents $\varphi$ near the current point $\mathbf{x}^\nu$ and $\Delta^\nu$ is the trust-region radius, that is adjusted at each iteration to produce a sufficiently decreasing approximation ($\varphi(\mathbf{x}^\nu+\mathbf{d}^\nu)<\varphi(\mathbf{x}^\nu)$).

The {\em reflective} algorithm uses as model function the quadratic function
\begin{equation}
m^\nu(\mathbf{d})=\frac{1}{2}\mathbf{d}^T\mathbf{B}^\nu \mathbf{d} + \mathbf{d}^T\nabla\varphi(\mathbf{x}^\nu),
\end{equation}
where $\mathbf{B}^\nu$ is the Hessian matrix $\nabla^2\varphi(\mathbf{x}^\nu)$ or an approximation to it \cite{nocedal2006numerical}. 
The minimization problem is restricted to $\mathbf{d}$ belonging to the two-dimensional subspace spanned by the gradient direction $\nabla \varphi(\mathbf{x}^\nu)$ and the Newton direction $\mathbf{B}^\nu \mathbf{d}=-\nabla \varphi(\mathbf{x}^\nu)$ \cite{byrd1988approximate}. For this algorithm, the Newton system is solved by applying the preconditioned conjugate gradient method \cite{branch1999subspace}.  

The trust-region {\em dogleg} algorithm, on the other hand, adopts \cite{conn2000trust}
\begin{equation}\label{eq:fsolve}
m^\nu(\mathbf{d})=\left\lVert \mathcal{H}(s^{\nu}) + \mathcal{H}'(s^{\nu})\mathbf{d}  \right\rVert_2^2.
\end{equation}
The update $\mathbf{d}$ is computed as a linear combination of the Cauchy and Newton steps, as presented in \cite{powell1968fortran}. 
The Cauchy step $\widetilde{\mathbf{d}}$ minimizes the model function $m^\nu$ in Eq. (\ref{eq:fsolve}) along the steepest descent direction. The Newton step is the unrestricted minimum of $m^{\nu}$ given by $\mathcal{H}'\widehat{\mathbf{d}}=-\mathcal{H}$ (i.e., the update defined in Eq. (\ref{BL2D_Newton8})). 
The dogleg algorithm chooses $\mathbf{d}=\widetilde{\mathbf{{d}}}+\chi(\widehat{\mathbf{d}}-\widetilde{\mathbf{d}})$, where $\chi$ is the largest value in $[0,1]$ such that $\parallel \mathbf{d}\parallel\leq\Delta^\nu$.

Our solver for the transport problem considers the three trust-region algorithms mentioned above. We perform comparisons of the approximations provided by them in the section with numerical results. Additionally, we consider for comparison the Newton method with a global under-relaxation factor of 0.5, in line with \cite{jenny2009unconditionally}. 
The Newton method with global under-relaxation is stable but requires significantly more iterations to converge when compared to the Newton method with the inflection-point strategy, where the under-relaxation is only applied locally.


\section{Numerical experiments}\label{sec_results}

In this section we present numerical experiments to study the performance of the sequential implicit solver using the MRCM for the approximation of two-phase flows. We compare the saturation approximations provided by the Newton method combined with the trust-region algorithms mentioned in subsection \ref{subsec:Trust-region algorithms}. 

The tolerance for the Newton step size is set to $\eta=10^{-6}$ in the $L^2$ norm, and the time is expressed in PVI (Pore Volume Injected) \cite{chen2006computational}. The relative permeabilities are given by $k_{ro}=(1-s)^2$ and $k_{rw}=s^2$, such that the fractional flow of water can be written as
\begin{equation}
f(s) = \dfrac{M s^2}{M s^2 + (1-s)^2},
\end{equation}
where $ M = {\mu_o}/{\mu_w}$. 
The numerical set-up in most of our simulations considers a flow established by imposing flux boundary conditions from left to right and no-flow at top and bottom. The domain is initially filled with oil, with water being injected at a constant rate. Source terms are zero and $M=10$. This is the configuration considered in the numerical studies unless stated otherwise.

We first compare the convergence of the saturation solution provided by the upwind method in the implicit and explicit versions. Then, we investigate the MRCM combined with the sequential implicit solver for different choices of permeability fields. 
Finally, we close our numerical experiments with an example that  considers gravity effects.

\subsection{Implicit versus explicit}\label{subsection:impli_expli}

\begin{figure}
		\centering
	\includegraphics[scale=0.38]{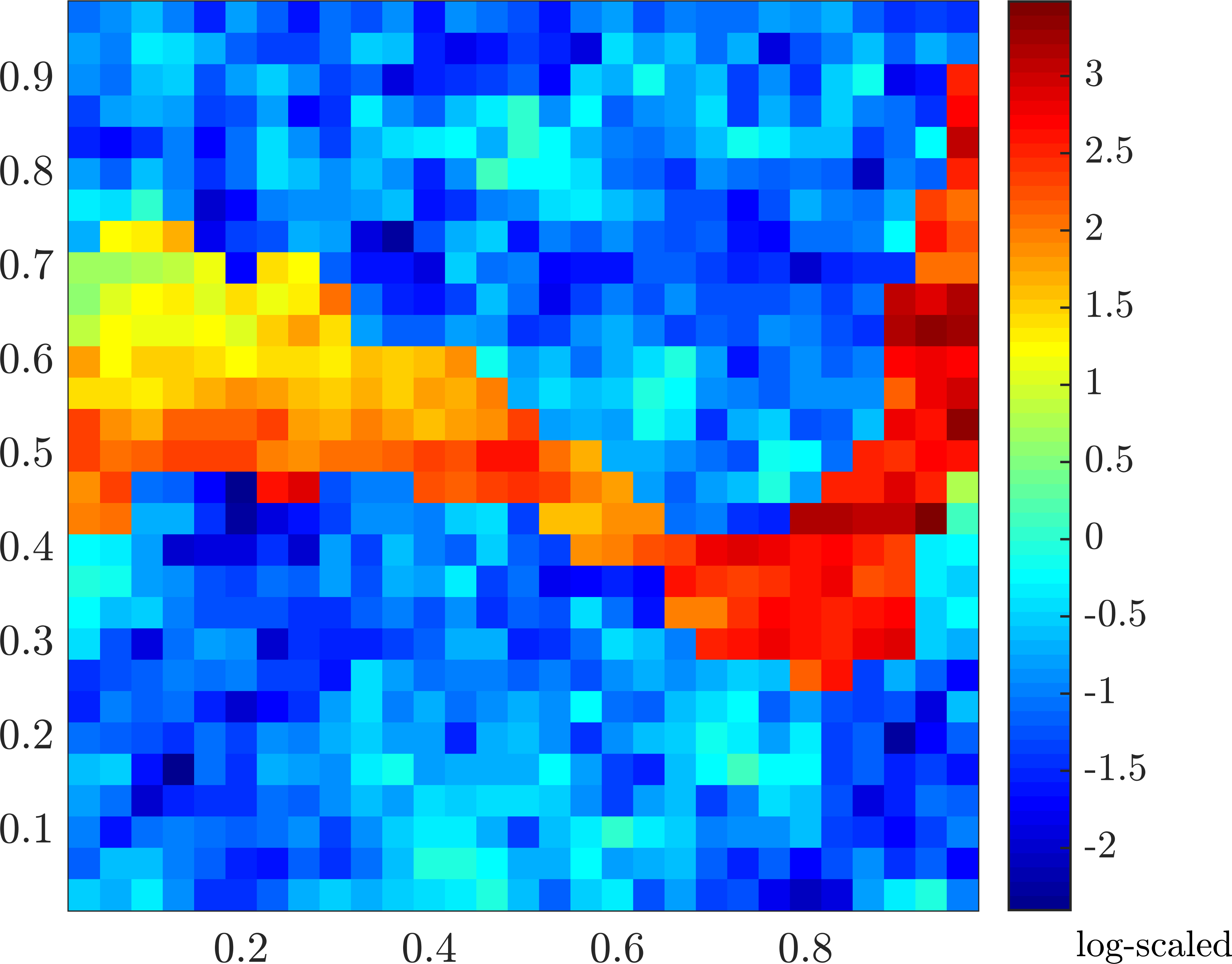}
	\includegraphics[scale=0.38]{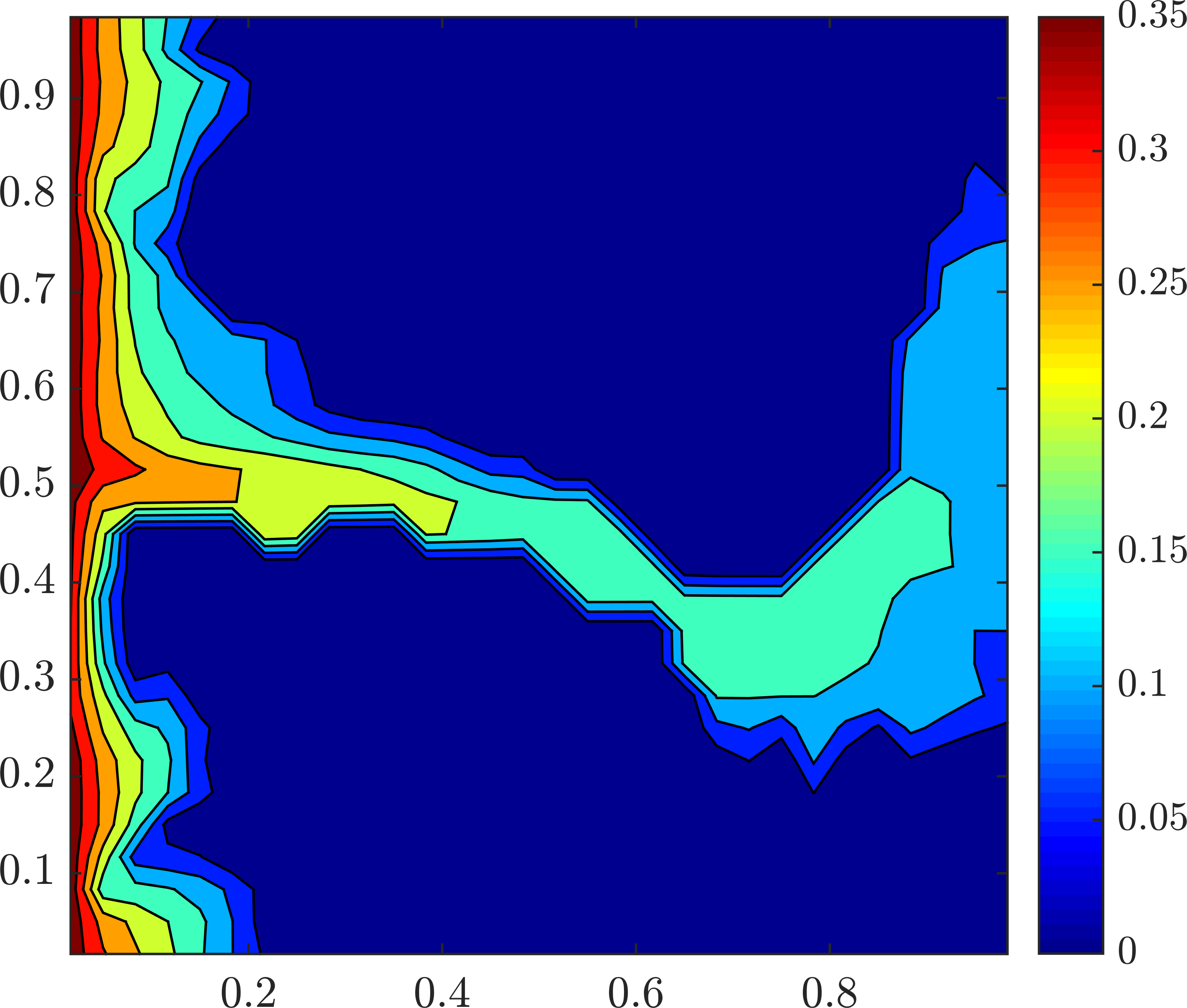}\quad
	\includegraphics[scale=0.55]{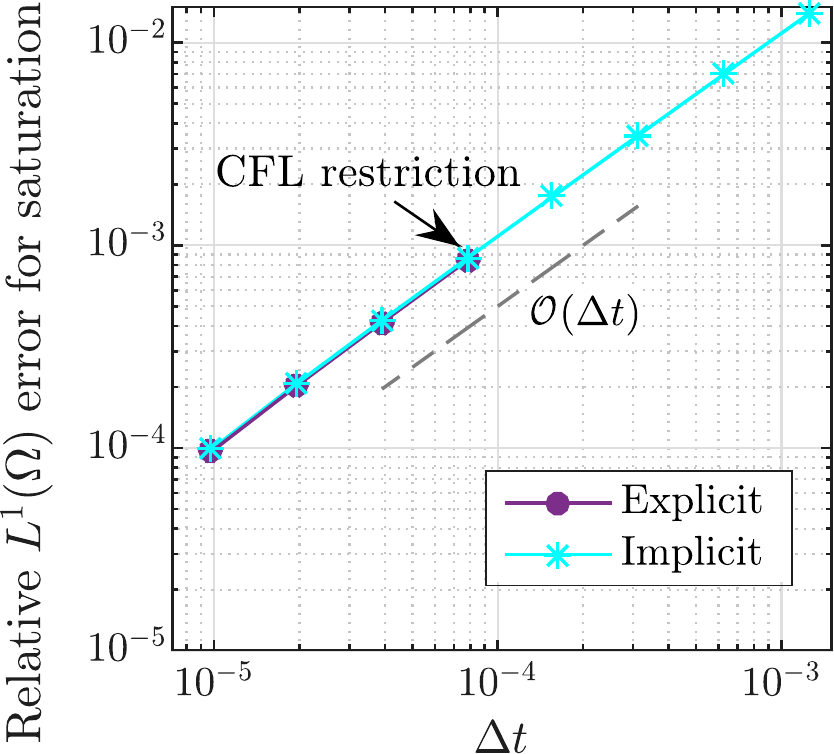}
	\caption{Log-scaled permeability field (left), saturation reference solution at the final time $T_{\textbf{PVI}}=0.0625$, and the relative $L^1(\Omega)$ error for saturation as a function of the time step size (right). The convergence behavior of the explicit and implicit schemes is the same for $\Delta t\leq \Delta t_{CFL}$ (linear slope), while only the implicit approximation is possible for larger sizes of time steps. }
	\label{fig:comp_expli_impli}
\end{figure}

The Newton method combined with the trust-region algorithms considered is unconditionally convergent, allowing for arbitrary sizes of time steps. Thus, the choice for the size of the time step is based only on accuracy requirements. 
One can take much larger time steps by using the implicit method instead of the explicit one for the transport equation in sequential approximations of two-phase flows. This is illustrated in Fig. \ref{fig:comp_expli_impli}, where we compare the convergence of the saturation solution provided by the upwind method in the implicit and explicit versions. For the implicit solution we use the Newton method with the inflection-point strategy and for the explicit case we consider the upwind method in an explicit operator splitting scheme (see \cite{douglas1997numerical} for additional discussion about the operator splitting framework). Here, the explicit approach fixes the same time step size for both elliptic and hyperbolic equations. 
 
 We consider a high-contrast permeability field in the domain $\Omega=[0,1]\times[0,1]$ (with $30\times30$ fine grid cells). We show in Fig. \ref{fig:comp_expli_impli} the log-scaled permeability filed (left), the saturation solution at the final time $T_{\textbf{PVI}}=0.0625$ (center), and the relative $L^1(\Omega)$ error for saturation as a function of the time step size (right). The time steps considered 
have size varying from $\Delta t= 9.75\times 10^{-6}$ to $\Delta t= 1.25\times10^{-3}$ (in PVI), while the reference solutions consider $\Delta t= 10^{-6}$ (in PVI, and satisfying the CFL condition). 
One can note that the errors are essentially the same for the time step sizes smaller than the CFL restriction ($\Delta t\leq \Delta t_{CFL}\approx 7.8\times 10^{-5}$). The explicit scheme cannot handle $\Delta t>\Delta t_{CFL}$, while the implicit method maintains the same behavior (linear slope) for the larger sizes of time steps. The average number of required Newton iterations per time step for the cases simulated vary from 2.36 for the smallest time step size to 6.5 for the largest one. Therefore, the implicit solver allows for approximating accurate solutions with large time step sizes, being a good choice to improve the efficiency of two-phase flow simulations.


\subsection{A Gaussian permeability field}\label{A comparison of implicit solutions}

In this example we consider a permeability field given by $K(\mathbf{x})=e^{4.5\xi(\mathbf{x})}$, where $\xi(\mathbf{x})$ is a self-similar Gaussian distribution having zero mean and covariance function given by $C(\mathbf{x},\mathbf{y})=|\mathbf{x}-\mathbf{y}|^{-1/2}$ \cite{glimm1993theory}. For this field, the permeability contrast is $K_{\max}/K_{\min} \approx 10^6$ and the computational grid has $64\times64$ cells in $\Omega=[0,1]\times[0,1]$.  
Concerning the MRCM, we consider linear interface spaces and set $\alpha(\mathbf{x})=1$. The domain decomposition considered has $4\times4$ subdomains, each one containing $16\times 16$ fine grid cells. In order to recover continuous fluxes at the interfaces of the skeleton we use the stitch downscaling procedure, we refer the reader to \cite{guiraldello2019downscaling} for details about this scheme.

\begin{figure}
\centering
\includegraphics[scale=0.8]{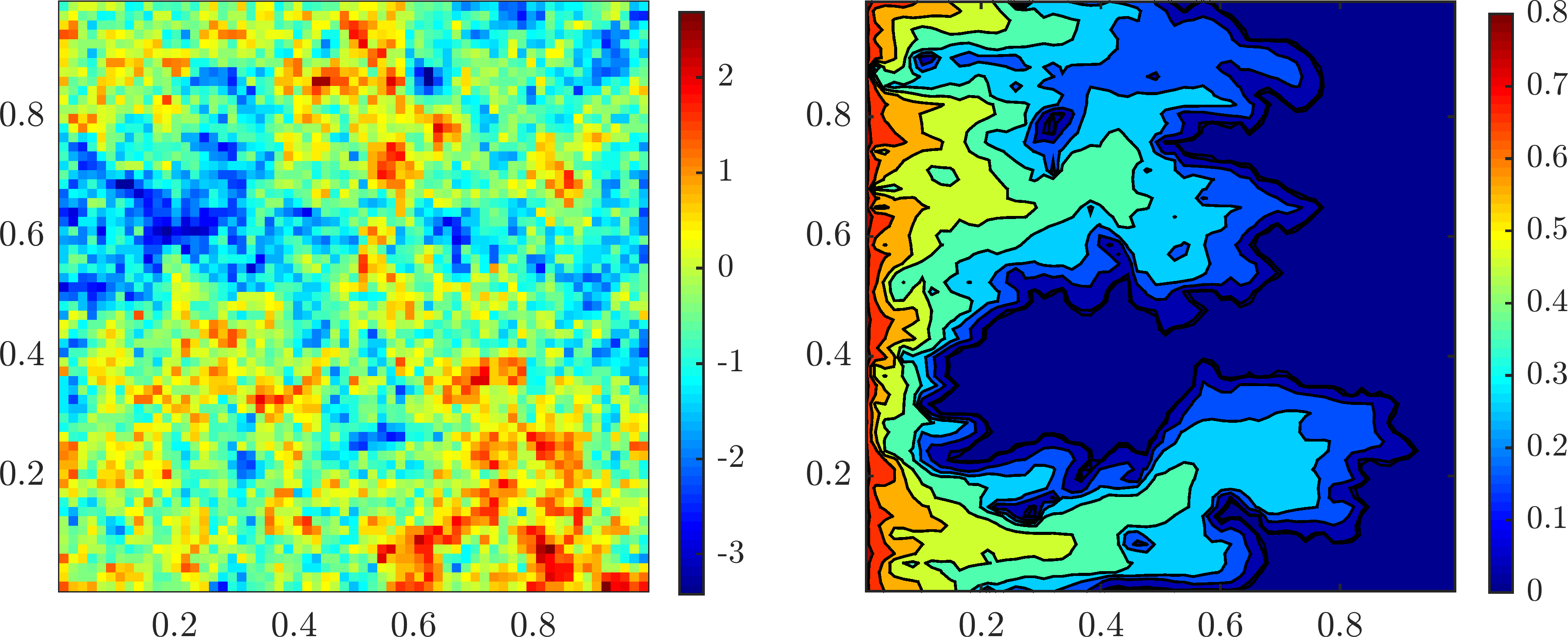}
	\caption{Log-scaled permeability field (left) and the saturation reference solution at the final time $T_{\text{PVI}}=0.2$ computed by the Newton method using the trust-region algorithm with the inflection-point strategy.}
	\label{fig:gaussiano}
\end{figure}

Figure \ref{fig:gaussiano} shows the log-scaled permeability field (left) and the saturation reference solution at the final time $T_{\text{PVI}}=0.2$ (right). This reference solution considers the Newton method using the trust-region algorithm with the inflection-point strategy for the transport problem. The time step for the reference solution is given by $\Delta t = 2 \times 10^{-5}\approx10^{-1}\Delta t_{CFL}$.

\begin{figure}
\centering
\includegraphics[scale=0.55]{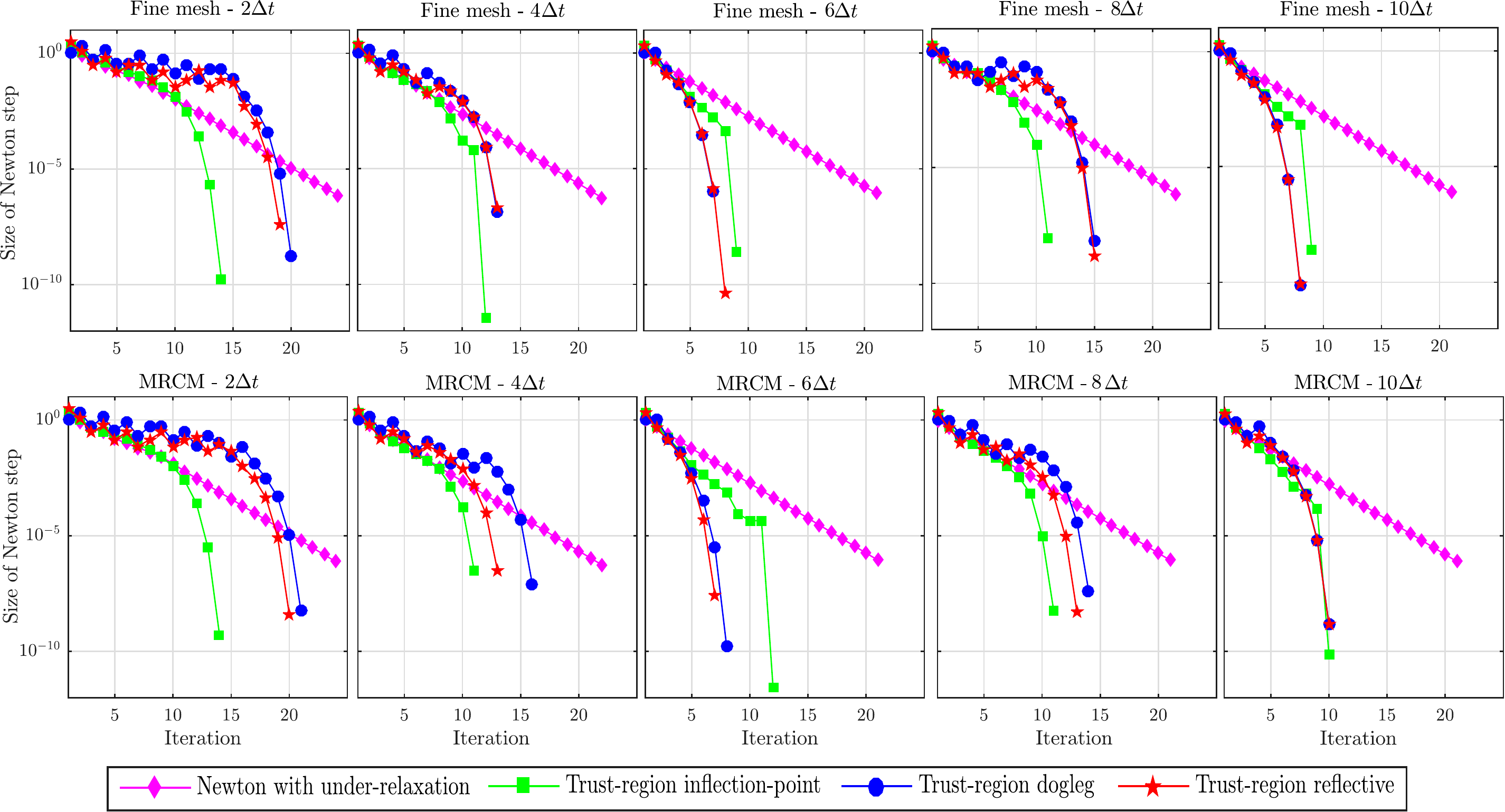}
	\caption{Convergence studies for one time step of size $\Delta t = 2\times10^{-2}\approx100\Delta t_{CFL}$. We show the size of the Newton step computed at times  $T_{\text{PVI}}=2\Delta t,\ 4\Delta t,\ 6 \Delta t,\ 8 \Delta t,\ 10 \Delta t$, considering the global under-relaxation and trust-region algorithms. The fine grid procedure (top row) and the MRCM (bottom row) have been considered to approximate the velocity field. The trust-region algorithms require significantly fewer iterations than the  under-relaxation technique.}
	\label{fig:gauss_step}
\end{figure}

The convergence results for the Newton method for a choice of $\Delta t = 2\times10^{-2} \approx100\Delta t_{CFL}$ (that corresponds to a total of 10 time steps) are shown in Fig. \ref{fig:gauss_step}. In this study, we compare the size of the Newton step provided by the trust-region algorithms. We also consider in the comparison the Newton method with global under-relaxation. The fine grid and MRCM procedures are considered for the elliptic updates. We show the size of the iterative Newton step computed at times $T_{\text{PVI}}=0.04,\ 0.08,\ 0.12,\ 0.16,\ 0.2$, i.e., $T_{\text{PVI}}=2\Delta t,\ 4\Delta t,\ 6 \Delta t,\ 8 \Delta t,\ 10 \Delta t$.
We note that the number of iterations needed by the trust-region algorithms are significantly smaller than the required by the  under-relaxation technique in all cases. The advantage of the the inflection-point strategy with respect to under-relaxation has been investigated in \cite{jenny2009unconditionally}. Here, we show that the trust-region dogleg and reflective algorithms can also be competitive to approximate the transport problem. Moreover, we find that by using the MRCM to update the velocity field we obtain similar number of Newton iterations.

\begin{figure}
\centering
\includegraphics[scale=0.6]{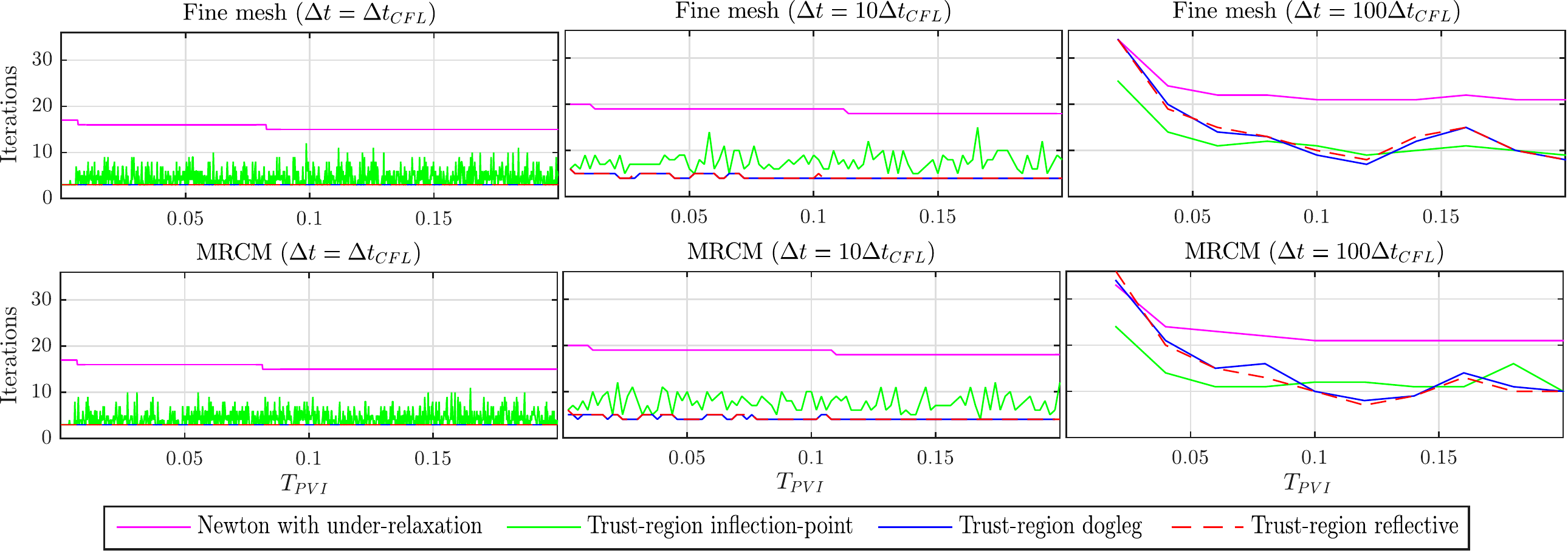}
	\caption{Number of Newton iteration required as a function of time (in PVI). Three time step choices as multiples of $\Delta t_{CFL}$ are shown. The fine grid procedure (top row) and the MRCM (bottom row) have been considered to approximate the velocity field. Note that the under-relaxation technique requires significantly more iterations than the trust-region algorithms.}
	\label{fig:iterations_gaussiano}
\end{figure}

\begin{figure}
\centering
\includegraphics[scale=0.75]{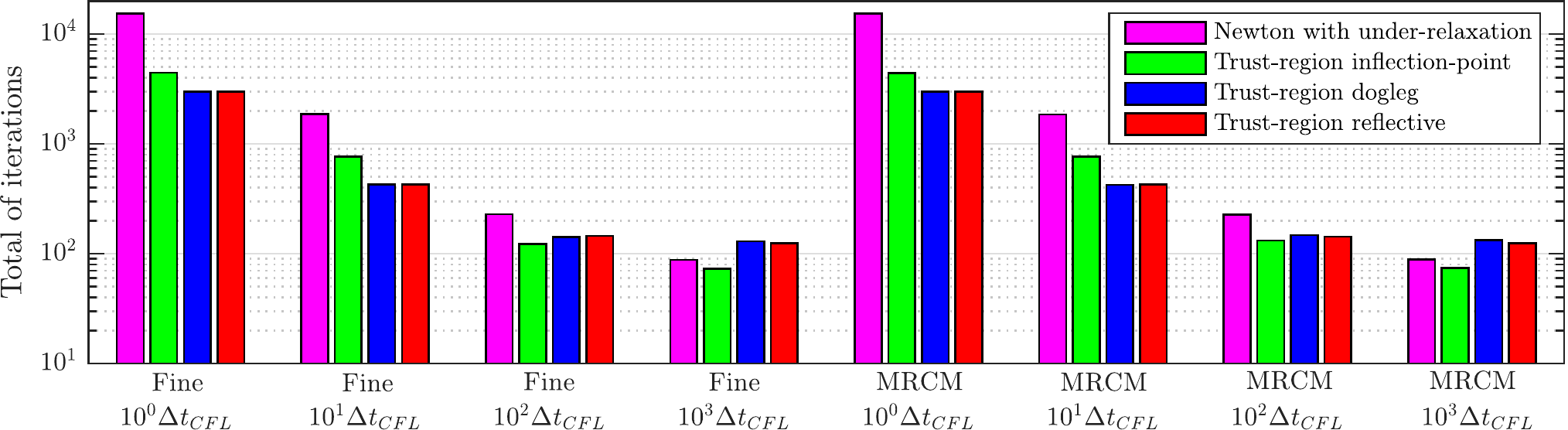}
	\caption{Total accumulated of iterations at the final time $T_{\textbf{PVI}}=0.2$ for different time step choices. The number  of iterations required by the procedures that use the fine grid and MRCM to compute the velocity field are shown, which are comparable. Note that the Newton method using the trust-region algorithm with the inflection-point strategy is the procedure that requires fewer iterations for the largest values of $\Delta t$ considered.}
	\label{fig:sum_iter_gaussiano}
\end{figure}

The behavior of the methods over time is shown in Fig. \ref{fig:iterations_gaussiano}, where we present the number of Newton iterations required in a simulation with three different time step choices. We present results for $\Delta t$ chosen as $\Delta t=\Delta t_{CFL}$, $\Delta t=10\Delta t_{CFL}$ and $\Delta t=100\Delta t_{CFL}$, that generate, respectively, a total of 1000, 100, and 10 time steps. In this case $\Delta t_{CFL}\approx 2\times 10^{-4}$. 
We note that the under-relaxation technique requires significantly more iterations than the trust-region algorithms. The trust-region dogleg and reflective algorithms require comparable numbers of iterations in all cases. For the values of $\Delta t=\Delta t_{CFL}$ and $\Delta t= 10\Delta t_{CFL}$, the trust-region dogleg and reflective algorithms require fewer iterations than the inflection-point strategy. For $\Delta t= 100\Delta t_{CFL}$, the performances of the three trust-region algorithms are similar, with a slight advantage for the inflection-point strategy. 
These results are summarized in Fig. \ref{fig:sum_iter_gaussiano}, where the total accumulated of Newton iterations is shown. 
We note an advantage in terms of the number of iterations for the trust-region algorithm with the inflection-point strategy when $\Delta t$ increases. 
 Both in Fig. \ref{fig:iterations_gaussiano} and Fig. \ref{fig:sum_iter_gaussiano} we note that the number of iterations required by the procedures that use the fine grid and MRCM to compute the velocity field is essentially the same. Therefore, the SI solver does not suffer from an increase in the number of iterations when combined with the MRCM.

\begin{figure}
\centering
\includegraphics[scale=0.75]{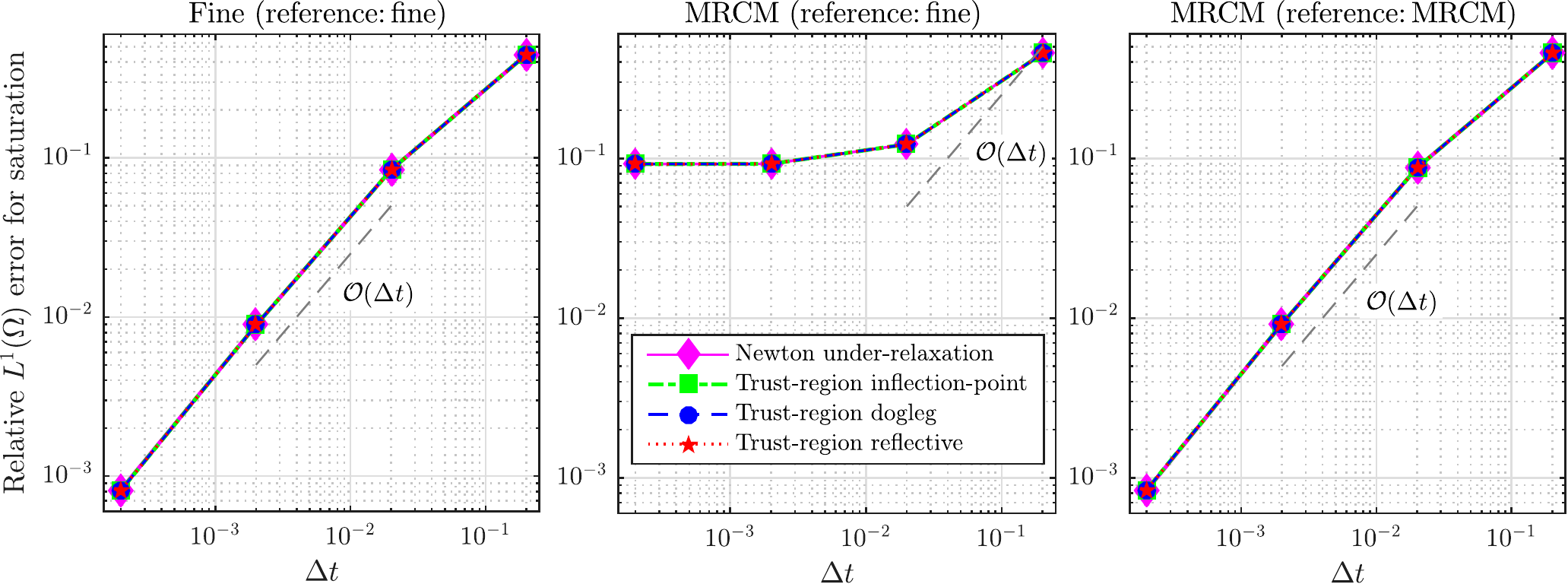}
\caption{Convergence in time by setting $\Delta t=10^{-1}\Delta t_{CFL}$ as the reference time step. Each method for the transport problem considers its corresponding reference solution. 
We show the convergence for the procedures that use the fine grid  (left) and MRCM (right) to compute the velocity field, considering as reference the fine grid velocity field and MRCM approximations, respectively. The convergence for the MRCM considering the fine grid velocity field as reference is shown at the center. We observe linear slope in all methods for the hyperbolic equation when excluding the multiscale errors of the velocity field.}
\label{fig:conv_gauss}
\end{figure}

We perform a convergence study by setting $\Delta t=0.1\Delta t_{CFL}$ as the reference time step.  
We estimate the errors of each method for the hyperbolic equation considering each reference solution computed with $\Delta t=0.1\Delta t_{CFL}$. In space, we consider the fine grid (with the fine grid velocity field approximation to be the reference) in Fig. \ref{fig:conv_gauss} (left), the MRCM (with the fine grid velocity field approximation as reference) in Fig. \ref{fig:conv_gauss} (center), and the MRCM (with the MRCM approximation as reference) in Fig. \ref{fig:conv_gauss} (right). 
 In the convergence study reported at the center of Fig. \ref{fig:conv_gauss}, the multiscale inaccuracies in the velocity field seem to be relevant when $\Delta t <2\times10^{-2}= 100\Delta t_{CFL}$, while the error of the transport process is dominant for the largest time step choices. We observe linear slope in all methods for the hyperbolic solver when excluding the multiscale errors, i.e., if only the time refinements are taking into account (the fine grid velocity field used in left of Fig. \ref{fig:conv_gauss} and MRCM at right of Fig. \ref{fig:conv_gauss}). We remark that the observed linear behavior is the expected slope, once we computed the time discretization by the first-order implicit Euler method.

\begin{figure}
\centering
\includegraphics[scale=0.75]{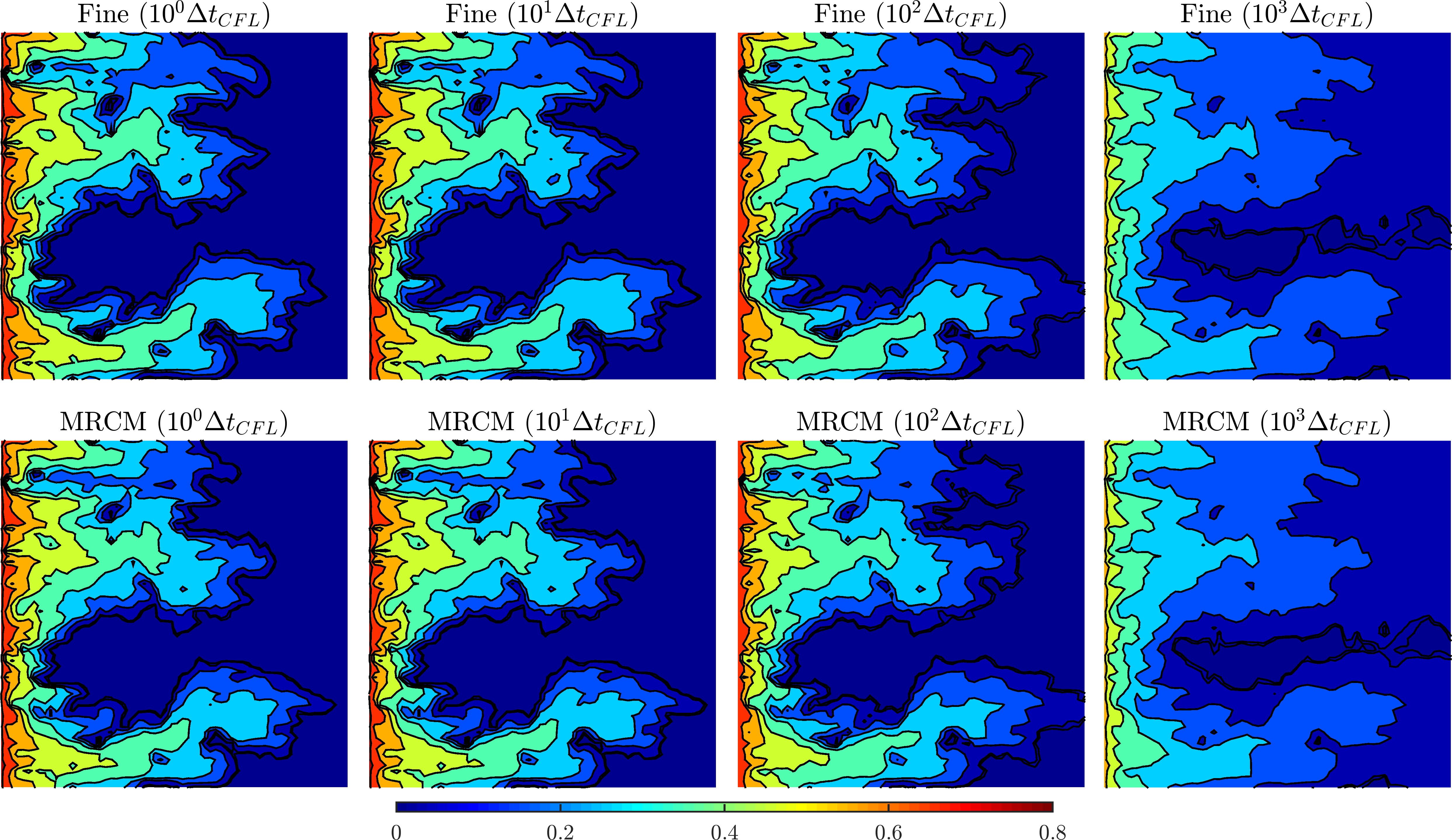}
\caption{Saturation profiles at time $T_{\textbf{PVI}}=0.2$ approximated by the Newton method using the trust-region algorithm with the inflection-point strategy. Different time step choices as multiples of $\Delta t_{CFL}$ are considered. We show the approximations that consider the fine grid (top row) and MRCM (bottom row) to compute the velocity field. For each choice of $\Delta t$, the MRCM and fine grid approximations are closely related.}
\label{fig:perfil_sat_gauss}
\end{figure}

To close this discussion, Fig.  \ref{fig:perfil_sat_gauss} shows a comparison of the saturation profiles at time $T_{\textbf{PVI}}=0.2$ approximated by the Newton method using the trust-region algorithm with the inflection-point strategy (since all the converged saturation solutions are the same). We show saturation maps obtained with different sizes of time steps. Note that all the approximations are consistent with the physics of the problem. In line with the convergence study we can conclude that the choice of the time step size is based only on accuracy requirements. For each choice of $\Delta t$, we note that the approximations that use the MRCM and fine grid velocity field are closely related. Note that there are no numerical artifacts originated strictly by the multiscale approximation of the velocity field.

\subsection{A channelized permeability field}\label{A multiscale solution}

\begin{figure}
\centering
\includegraphics[scale=0.42]{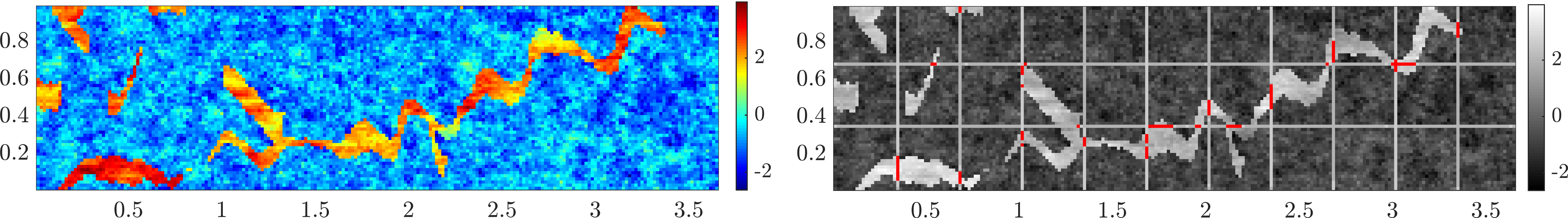}
	\caption{Permeability field (log-scaled) from layer number 36 of the SPE10 project (left) and a map of the absolute permeability variations at the boundaries of the subdomains (right). The red color identifies the high-permeability channel, where $\alpha(\mathbf{x})=10^{-2}$ is set and the physics-based interface spaces for pressure are defined.}
	\label{spe10}
\end{figure}

The next example considers the layer number 36 of the SPE10 project \cite{christie2001tenth}, that has a highly-permeable channel and permeability contrast of $K_{\max}/K_{\min} \approx 10^6$, see Fig. \ref{spe10} (left). The domain for this example is $\Omega=[0,11/3]\times[0,1]$ with $220\times60$ fine grid cells. 
For this high-contrast channelized formation, we apply the physics-based interface space for pressure to better represent the solution in the high-permeability channel. 
We investigate the accuracy of the MRCM combined with the SI approach. 
The domain decomposition considered contains $11\times 3$ subdomains with $20\times 20$ cells in each one of them. The flux interface space is linear, as well as the pressure space at the interfaces that do not cross the high-permeability channel. We use the adaptive version of the MRCM \cite{bifasico} by setting $\alpha(\mathbf{x})=10^{-2}$ at the interfaces that cross the high-permeability channel and $\alpha(\mathbf{x})=10^{2}$ at the remaining interfaces. Figure \ref{spe10} (right) shows a map of the absolute permeability variations at the boundaries of the subdomains. The red color identifies the high-permeability channel, where $\alpha(\mathbf{x})=10^{-2}$ is set and the physics-based interface spaces for pressure are defined.

\begin{figure}
\centering
\includegraphics[scale=0.75]{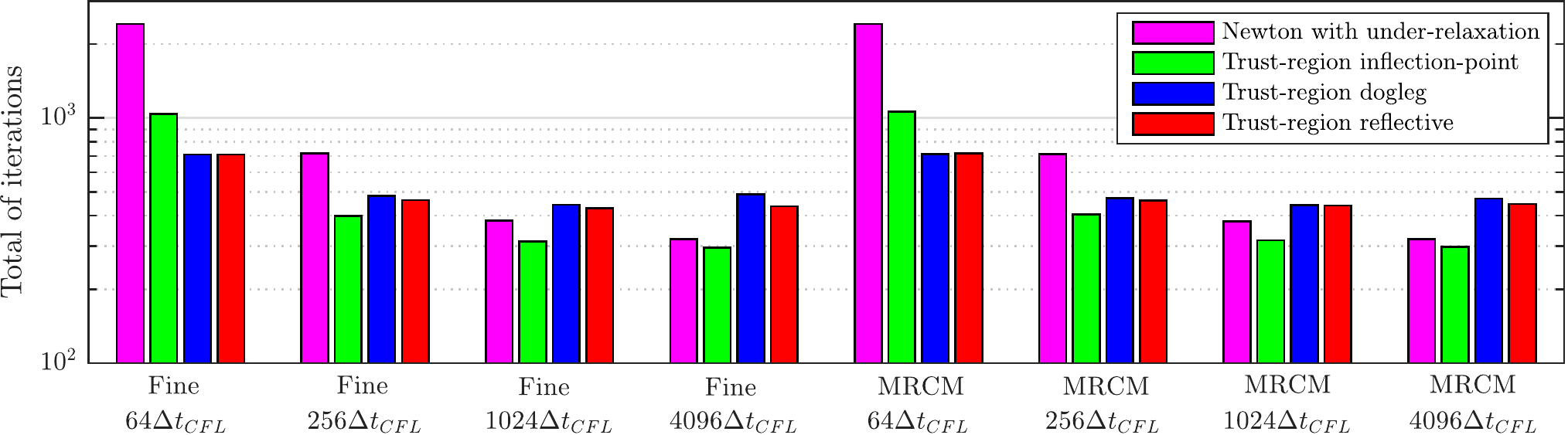}
	\caption{Total accumulated of Newton iterations until time $T_{\textbf{PVI}}=0.11$ for the SI scheme combined with the fine grid solution and MRCM. Different sizes of time steps are considered. The number of iterations required by the MRCM is comparable to the needed by the fine grid solution in all cases.}
	\label{fig:sum_iter_36_mrc}
\end{figure}

Figure \ref{fig:sum_iter_36_mrc} shows the total accumulated of Newton iterations until time $T_{\textbf{PVI}}=0.11$ for different sizes of time step taken as multiples of $\Delta t_{CFL}\approx 1.3\times 10^{-5}$. We start with $\Delta t=64\Delta t_{CFL}$ (that corresponds to a total of 128 time steps) and multiply by four until $\Delta t=4096\Delta t_{CFL}$ (that corresponds to a total of 2 time steps). We show results for the SI scheme combined with the fine grid velocity field and MRCM. 
We note a clear advantage in the number of Newton iterations for the trust-region algorithm with the inflection-point strategy when $\Delta t$ increases.
The trust-region reflective and dogleg schemes are more competitive for sizes of time step chosen of the order of $10\Delta t_{CFL}$, while the inflection-point strategy is the best choice for sizes of time step of the order of $100\Delta t_{CFL}$ or $1000\Delta t_{CFL}$. The under-relaxation technique, as expected, requires more Newton iterations to converge than the trust-region algorithm with the inflection-point strategy. However, it performs better than the trust-region reflective and dogleg algorithms for large sizes of time step. The number of iterations required by the procedure that uses the MRCM to compute the velocity field is comparable to that needed by the fine grid solution in all cases.

\begin{figure}
\centering
\includegraphics[scale=1.1]{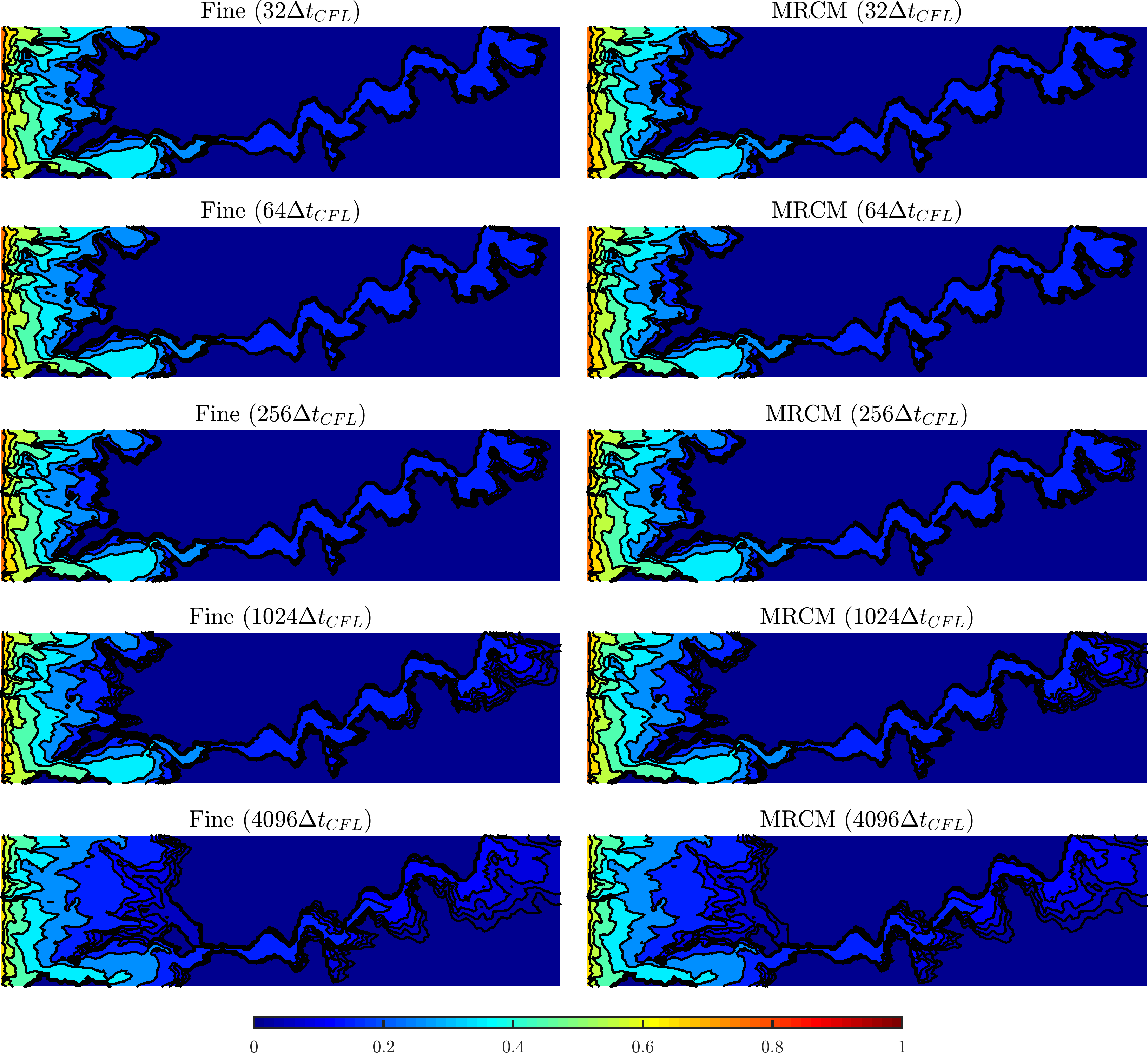}
\caption{Saturation profiles at time $T_{\textbf{PVI}}=0.11$ approximated by the procedures that use the fine grid and MRCM to compute the velocity field, considering different sizes of time step. For each choice of $\Delta t$ the approximations provided by the MRCM and fine grid procedure are closely related.}
\label{fig:sat_perfil_36_mrc}
\end{figure}

A comparison of the saturation profiles at time $T_{\textbf{PVI}}=0.11$ approximated by the fine grid procedure and MRCM is displayed in Fig. \ref{fig:sat_perfil_36_mrc}. The transport problem considers the Newton method using the trust-region algorithm with the inflection-point strategy. We remark that after the convergence of the Newton method all the algorithms considered (the under-relaxation technique, inflection-point, trust-region dogleg, and trust-region reflective) provide similar accuracy. We present saturation maps given by different sizes of time steps. For each choice of $\Delta t$, we note that the approximations that use the fine grid velocity field and MRCM are closely related.

\begin{figure}
\centering
\includegraphics[scale=0.65]{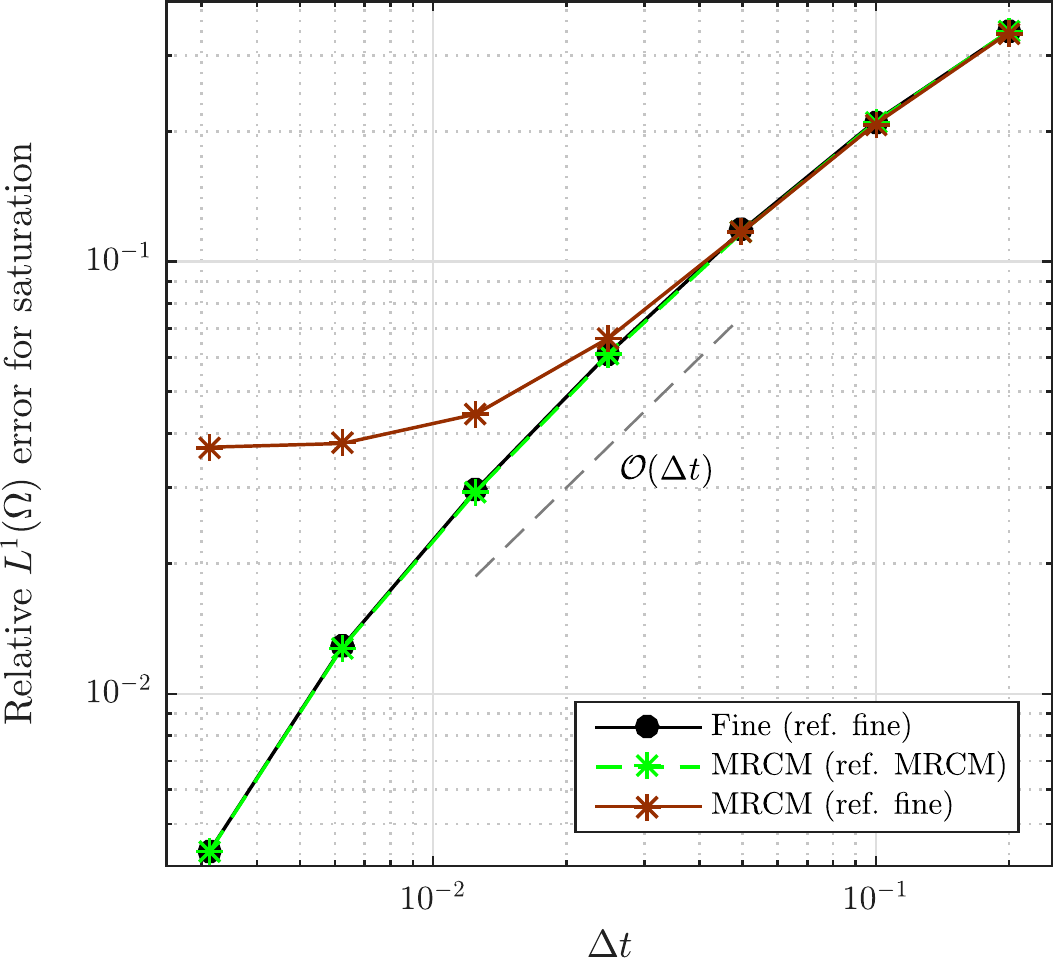}
\caption{Convergence in time by setting the solution computed with $\Delta t=32\Delta t_{CFL}$ as reference. We compare the errors of the procedures that use the fine grid and MRCM to compute the velocity field, the latter considering its corresponding reference (in space) and the fine grid reference solution. The multiscale inaccuracies of the velocity field are relevant for choices of $\Delta t <2.5\times 10^{-2}$ (that corresponds to $\Delta t <512\Delta t_{CFL}$), while the error of the transport procedure is dominant for larger time step choices.}
\label{fig:conv_Newton_36_mrc}
\end{figure}

The study of the convergence in time is reported in Fig. \ref{fig:conv_Newton_36_mrc}, where we consider as reference the approximation computed with $\Delta t =32\Delta t_{CFL}$. Here, we choose the Newton method using the trust-region algorithm with the inflection-point strategy to report the convergence in time since all the hyperbolic solvers presented the same behavior. Linear slope is attained when the error of each procedure (fine grid velocity field or MRCM) considers its corresponding spatial reference solution. The MRCM errors with respect to the fine grid solution show that the multiscale inaccuracies of the velocity field are relevant when $\Delta t <2.5\times 10^{-2}$ (that corresponds to $\Delta t <512\Delta t_{CFL}$), while the error of the transport procedure is dominant for larger time step choices. If we compute only the  multiscale error, i.e., the error of the MRCM by considering as reference the respective fine grid velocity field with the same time discretization, we obtain an error of the order $3\times 10^{-2}$ for all $\Delta t$, which is consistent with the dominant multiscale error observed on the convergence curve.

\begin{figure}
\centering
\includegraphics[scale=0.7]{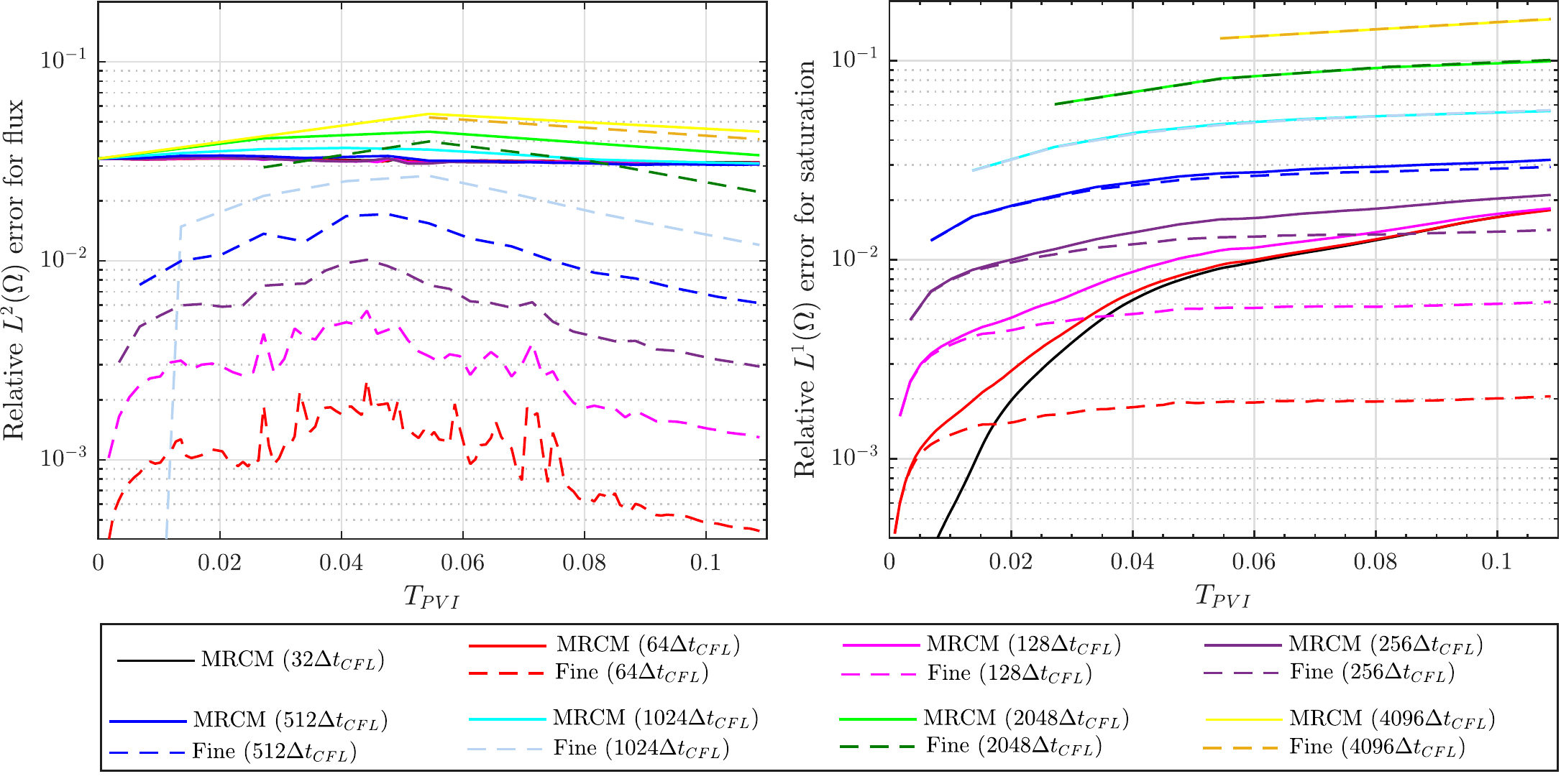}
\caption{Relative errors of flux (left) and saturation (right) for the procedures that use the MRCM (solid lines) and fine grid (dashed lines) to compute the velocity field, considering different time step choices. The reference solution is the fine grid approximation with $\Delta t=32\Delta t_{CFL}$. We note that the error curves do not vary significantly over time.}
\label{fig:Newton_mrc_erro_time}
\end{figure}

Finally, we report in Fig. \ref{fig:Newton_mrc_erro_time} the errors of the previous experiment as a function of time. We show relative errors of flux and saturation for the procedures that use the MRCM (solid lines) and fine grid approximations (dashed lines) to compute the velocity field, considering different time step choices. The reference is the fine grid solution with $\Delta t=32\Delta t_{CFL}$. 
The flux error is computed as usual: at each time we divide the $L^2$ norm of the difference by the $L^2$ norm of the reference at the same time. The saturation error (for this plot) divides the $L^1$ norm of the differences by the maximum absolute of the reference on time, avoiding divisions by very small values at the beginning of the simulation.
We note that the error curves do not vary significantly over time. 
This result shows that the observations from the convergence study for saturation (at time $T_{\textbf{PVI}}=0.11$) are maintained throughout the simulation. 
The flux errors of the MRCM are essentially the same for all $\Delta t$, whereas the fine grid errors decrease with the size of the time step.

\subsection{A homogeneous medium with fingering instability}\label{Glimm}

In this subsection we test the sequential implicit solver for a challenging problem with a fingering instability in a homogeneous medium. In contrast to the previous examples, here we do not consider heterogeneity, however the high nonlinearity of the coupling of flow and transport generates an unstable oil-water interface.

\begin{figure}
\centering
\includegraphics[scale=0.9]{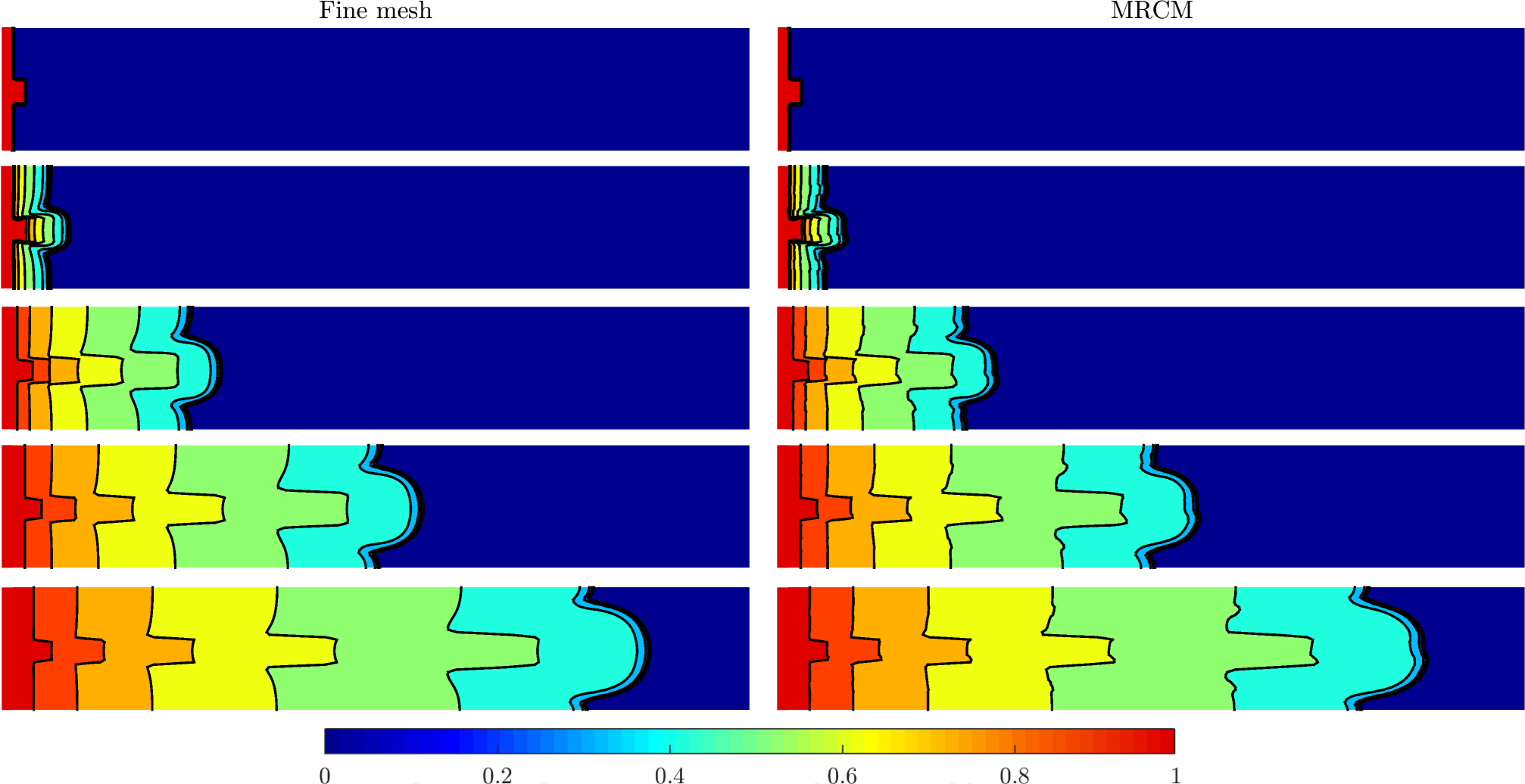}
\caption{Saturation solution computed by the fine grid procedure (left) and the MRCM (right), both combined with the SI solver. We show the profiles at times $T_{\text{PVI}}=0.00,\ 0.03,\ 0.15,\ 0.34,\ 0.60$, from top to bottom. Here, $\Delta t=\Delta t _{CFL}$. Note that the finger growing is well captured by both methods.}
\label{fig:Glimm_1}
\end{figure}

The numerical set-up considers a slab geometry with flow established from left to right by imposing pressure $p = 0$ on the left and $p = -10^{4}$ on the right along with no-flow at top and bottom. The domain $\Omega=[0,3]\times[0,1/2]$, with $300\times50$ fine grid cells, has an initial water front at left, while the rest of the reservoir is filled with oil. The water front contains a small perturbation at the center, as shown in Fig. \ref{fig:Glimm_1}. No source terms are considered. 
This geometry generates a finger that evolves in time, characterizing a 2D Riemann problem \cite{glimm1981numerical}. 
The nonlinearity, and hence, the physical instabilities of this problem are connected to the viscosity ratio value. Here, we choose $M=4$ in line with \cite{glimm1981numerical} and \cite{furtado2003crossover}, where the authors have shown that the critical value for unstable flows is $M \approx 2.657$.   
The MRCM approximation for this homogeneous porous medium uses linear interface spaces and set $\alpha(\mathbf{x})=1$. A domain decomposition of $15\times 5$ subdomains with $20\times 10$ cells into each one is considered.

Figure \ref{fig:Glimm_1} shows the saturation approximations computed by the SI solver combined with the fine grid (left) and MRCM (right). The profiles at times $T_{\text{PVI}}=0.00,\ 0.03,\ 0.15,\ 0.34,\ 0.60$, from top to bottom are shown. For the current example, we approximate the transport problem by the Newton method using the trust-region algorithm with the inflection-point strategy.  
The saturation illustrated in Fig. \ref{fig:Glimm_1} is our reference solution, computed with $\Delta t=\Delta t _{CFL}\approx 5.82\times10^{-6}$. Although some inaccuracies appear in the MRCM approximation when compared to the fine grid solution, the finger growth is well captured by both methods.

\begin{figure}
\centering
\includegraphics[scale=1]{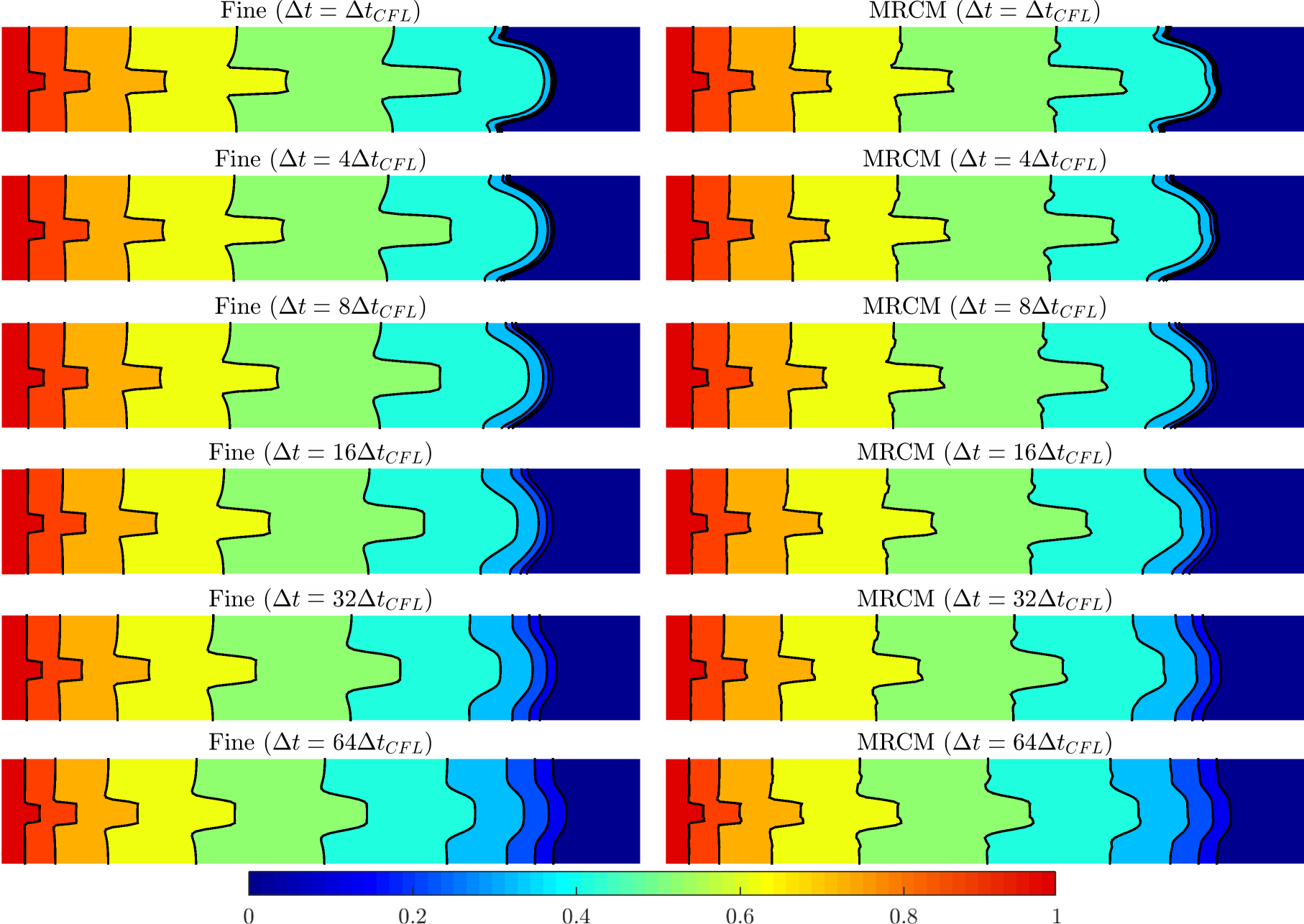}
\caption{Saturation profiles at time $T_{\textbf{PVI}}=0.6$ approximated by the procedures that use the fine grid (left) and MRCM (right) to compute the velocity field, considering different sizes of time step. Relevant inaccuracies in the water front appear for values of $\Delta t$ as of the order of $10\Delta t_{CFL}$.}
\label{fig:Glimm_2}
\end{figure}

A comparison of the saturation profiles at time $T_{\textbf{PVI}}=0.6$ approximated by the procedures that use the fine grid and MRCM to compute the velocity field is shown in Fig. \ref{fig:Glimm_2}, where different sizes of time steps are considered. For each choice of $\Delta t$, the fine grid and MRCM approximations are close, the latter presenting only small inaccuracies. Note that this problem presents relevant inaccuracies related to the transport approximations for choices of $\Delta t$ of the order of $10\Delta t_{CFL}$. 

Although the Newton method in the SI solver is unconditionally convergent, the accuracy for this problem with fingering instability drops quickly when $\Delta t$ increases. We intend to choose large sizes of time step, aiming at computational efficiency. Therefore, we combine the MRCM with the Sequential Fully-Implicit (SFI) scheme \cite{jenny2006adaptive}, which is adequate to simulate more complex models in which the semi-implicit treatment of the velocity can generate low accuracy \cite{aziz1979petroleum}.

\subsubsection{Sequential fully-implicit approximation}

The Sequential Fully-Implicit \cite{jenny2006adaptive} method consists of an outer loop to solve the coupled problems of flow and transport at each time step and the inner (nonlinear) Newton loop to solve the implicit transport problem.

Concerning the SI solver explained in section \ref{sec_SFI}, an additional external loop is added to advance from time $t^n$ to time $t^{n+1}$, where a sequence of updates of the two equations is executed: the elliptic equation for pressure and flux Eq. (\ref{Darcy}), and the transport equation for saturation Eq. (\ref{BL2D}) (by using the Newton method). The sequential updates of both equations are repeated until the maximum absolute change in the saturation between successive iterations is less than a tolerance criterion. In our experiments, we set the tolerance for the external loop as $10^{-4}$, while the tolerance for the inner Newton loop is the same previously considered $\eta=10^{-6}$. 
We remark that if a single iteration of the external loop is performed, the SI algorithm is recovered.

\begin{figure}
\centering
\includegraphics[scale=1]{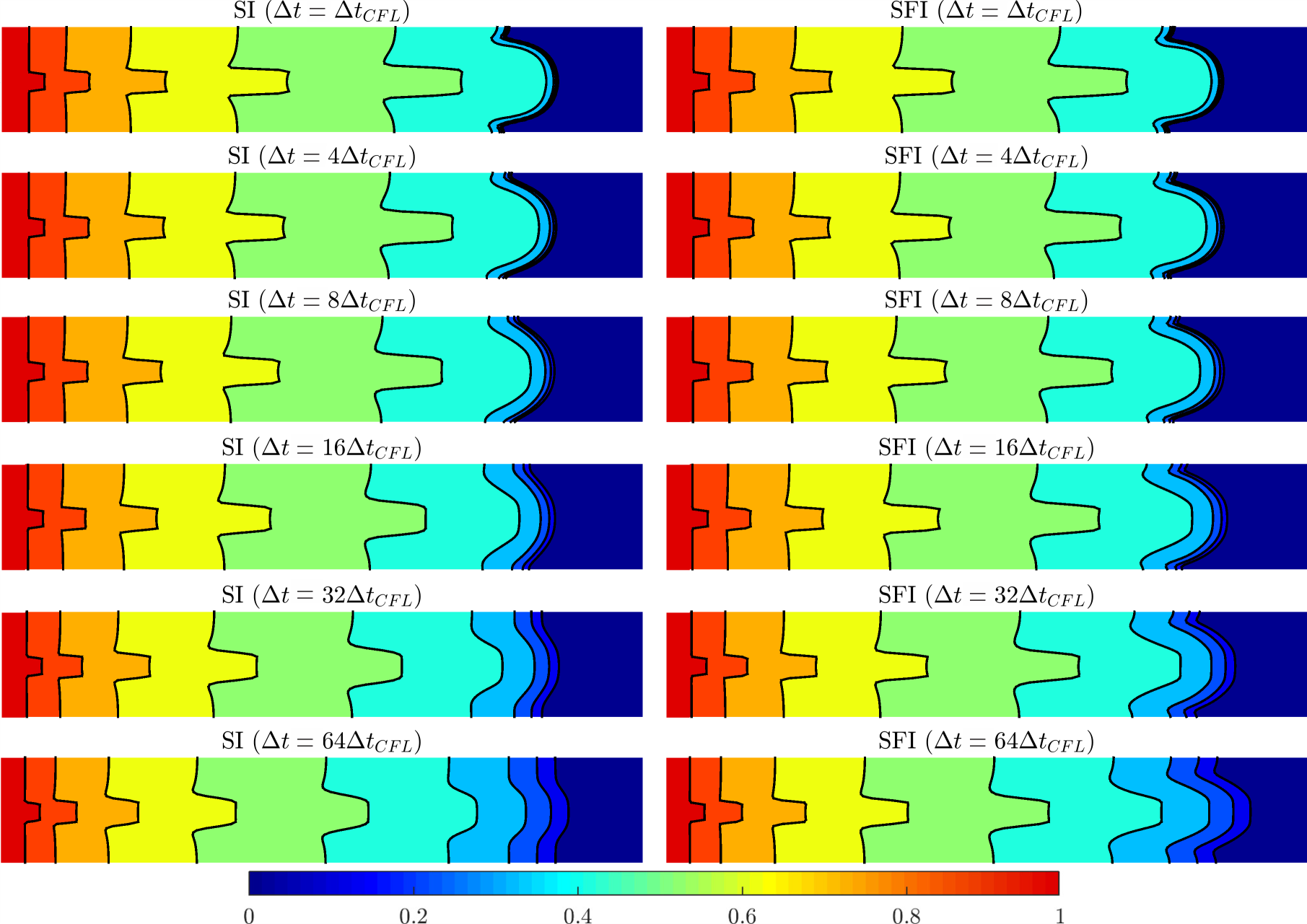}
\caption{Saturation profiles at time $T_{\textbf{PVI}}=0.6$ approximated by the fine grid velocity field combined with the SI (left) and SFI (right) schemes. Different sizes of time steps are considered. The SFI approximations are more accurate than the SI ones.}
\label{fig:Glimm_3}
\end{figure}

\begin{figure}
\centering
\includegraphics[scale=1]{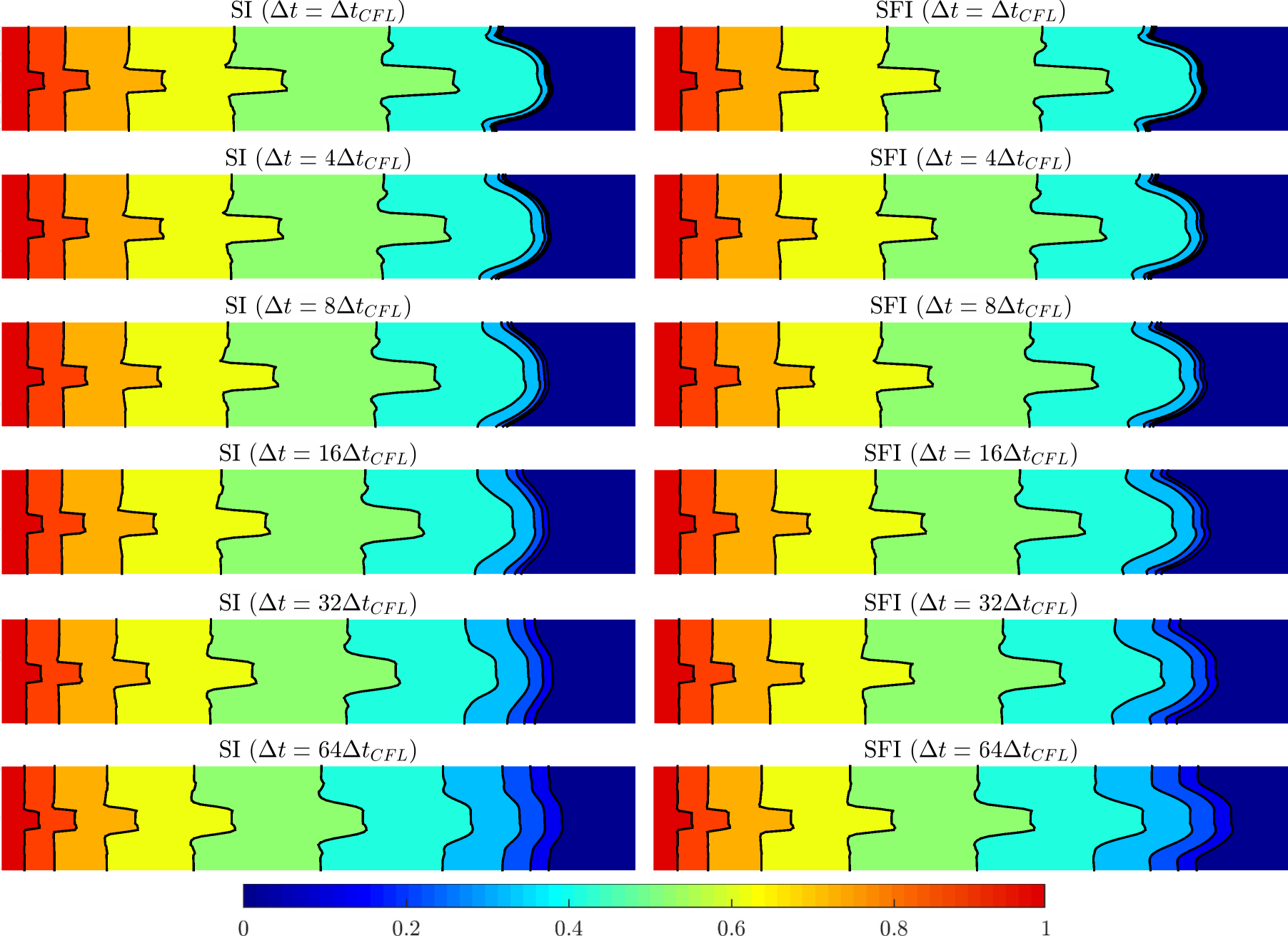}
\caption{Saturation profiles at time $T_{\textbf{PVI}}=0.6$ approximated by the MRCM combined with the SI (left) and SFI (right) schemes. Different sizes of time steps are shown. The MRCM works properly when combined with the SFI scheme.}
\label{fig:Glimm_4}
\end{figure}

Figure \ref{fig:Glimm_3} shows a comparison of the saturation profiles at time $T_{\textbf{PVI}}=0.6$ approximated by the SI and SFI schemes for different sizes of time steps. In this figure, we show results for the procedure that uses the fine grid velocity field. The reference solutions for each case are obtained by taking $\Delta t=\Delta t_{CFL}$. Note that the reference solutions provided by the SI and SFI solvers are similar. For each choice of $\Delta t>\Delta t_{CFL}$, the SFI approximations are more accurate than the SI ones. The SFI scheme presents inaccuracies for choices of $\Delta t \gtrsim 10\Delta t_{CFL}$ as well as the SI scheme, however the oil-water interface is better captured by the SFI solver. The corresponding comparison between the SI and SFI schemes combined with the MRCM is shown in Fig. \ref{fig:Glimm_4}. We note that the MRCM works properly when combined with the SFI scheme. The same observations about the relation between the SI and SFI approximations can be drawn if we combine them with the fine grid solution or MRCM.

We show a convergence study for the SI and SFI schemes by setting their respective solutions computed with $\Delta t=\Delta t_{CFL}$ as references. The fine grid and MRCM are used to compute the velocity field, the latter considering its corresponding reference (in space) and the fine grid reference solution. Essentially the same behavior (linear slope) is observed in all the curves. The SFI solver produces errors slightly lower than the SI scheme. However, the differences between the SI and SFI solutions are better observed in the saturation maps, where we note more detailed information when compared to the global error norms.

\begin{figure}
\centering
\includegraphics[scale=0.75]{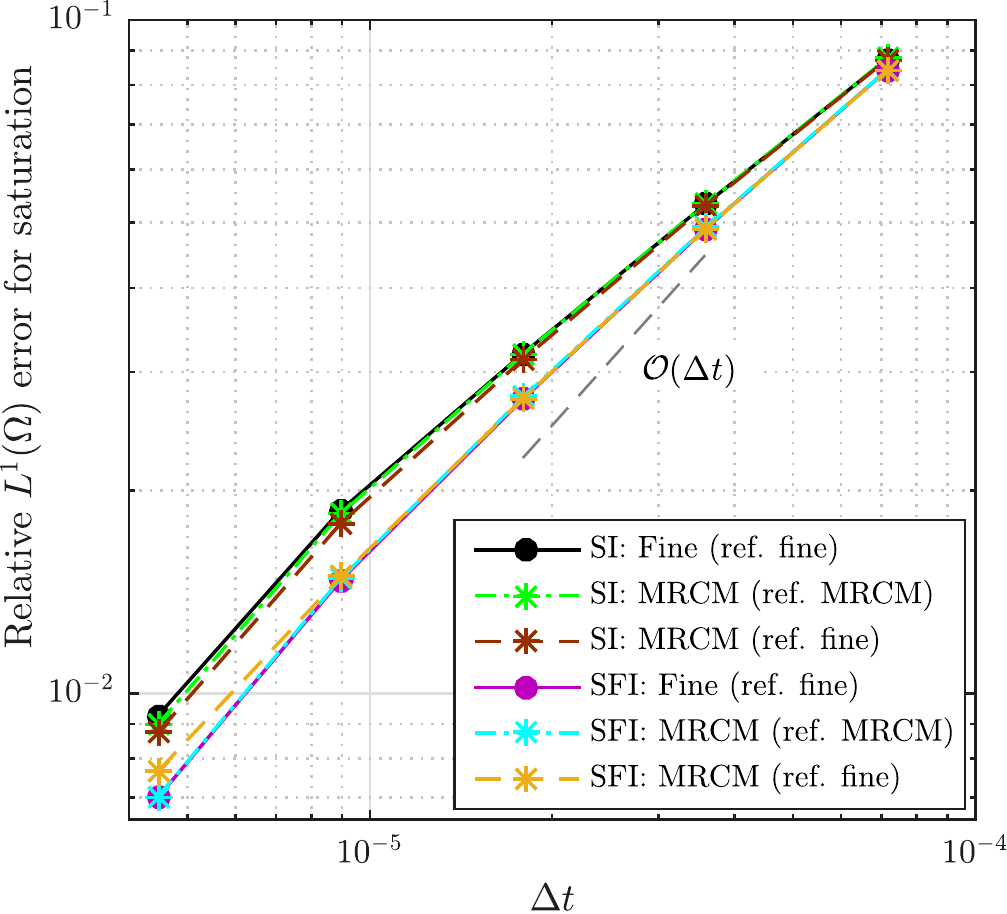}
\caption{Convergence in time of the SI and SFI schemes by setting the solution computed with $\Delta t=\Delta t_{CFL}$ as reference. The fine grid and MRCM are used to compute the velocity field, the latter considering its corresponding reference (in space) and the fine grid reference solution. All the curves present linear slope.}
\label{fig:Glimm_5}
\end{figure}

\begin{figure}
\centering
\includegraphics[scale=0.7]{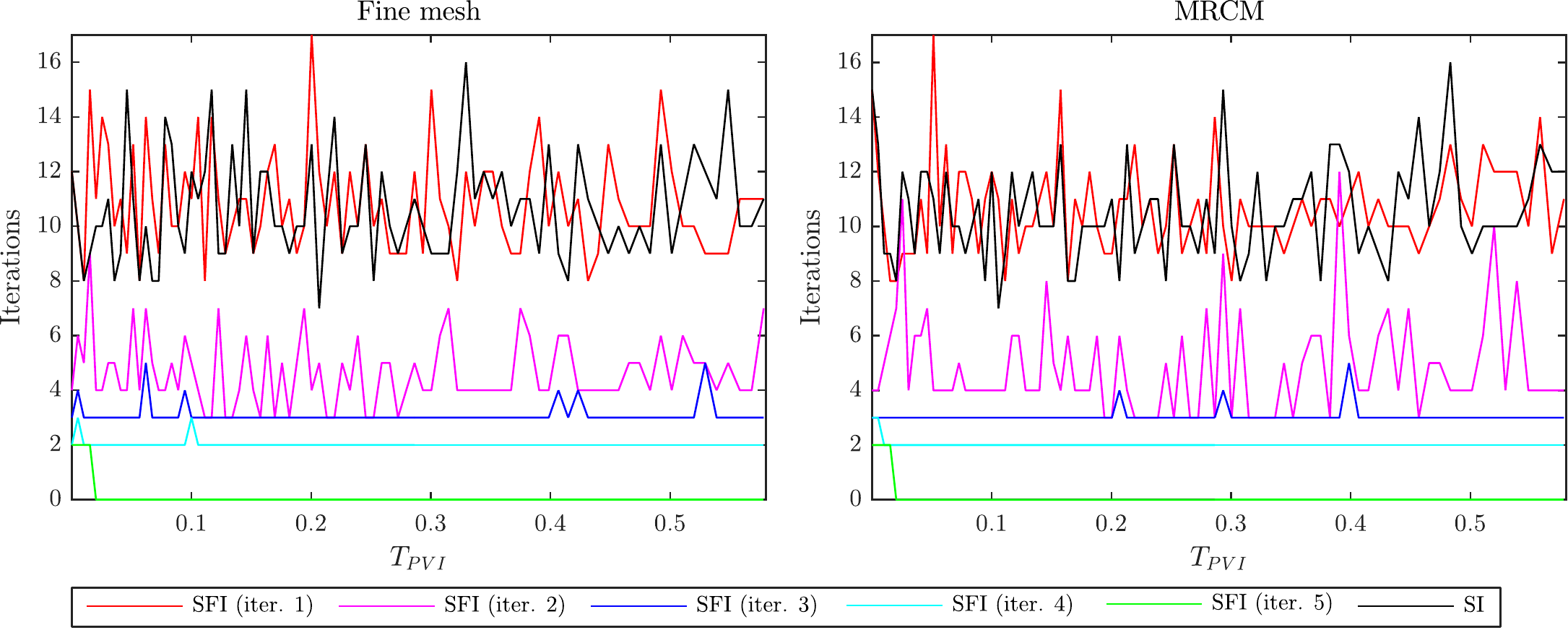}
\caption{Number of iterations required by the Newton method as a function of time (in PVI) for $\Delta t= 16\Delta t_{CFL}$. The numbers for the SI and SFI schemes are presented, the latter being composed of five external iterations. The fine grid procedure (left) and the MRCM (right) have been considered to approximate the velocity field. Note that the first external iteration of the SFI requires a similar number of Newton iterations to the SI scheme.}
\label{fig:Glimm_6}
\end{figure}

From the saturation maps we can conclude that the SFI solver is more accurate, and captures well the fingering instabilities when compared to the SI scheme. In terms of computational efficiency, the SFI is clearly more expensive than the SI (see \cite{ali2020multiscale, mpm2p, rochaenhanced} that would significantly decrease the cost of the SFI procedure). However, larger time steps can be chosen when using the SFI solver. We show in Fig. \ref{fig:Glimm_6} the number of iterations required by the Newton method to approximate the 2D Riemann problem for the intermediate size of time step $\Delta t= 16\Delta t_{CFL}$. The number of Newton iterations for the SI and SFI schemes are presented, the latter being composed of five external iterations. Note that the first external iteration of the SFI requires a similar number of Newton iterations to the SI scheme. The subsequent SFI iterations require a smaller number of Newton iterations: the second requires around 5; the third requires around 3; the fourth requires around 2; the fifth is only necessary at the beginning of the simulation. These observations are essentially the same for both procedures in the case of the fine grid velocity field Fig. \ref{fig:Glimm_6} (left) and MRCM Fig. \ref{fig:Glimm_6} (right). Thus, to attain better accuracy than that provided by the SI scheme, some external iterations in the SFI solver are necessary. For this 2D Riemann problem, a maximum of four external iterations was required when $\Delta t= \Delta t_{CFL}$, while a maximum of seven external iterations was required in the case of $\Delta t= 16\Delta t_{CFL}$. 

\begin{figure}
\centering
\includegraphics[scale=0.8]{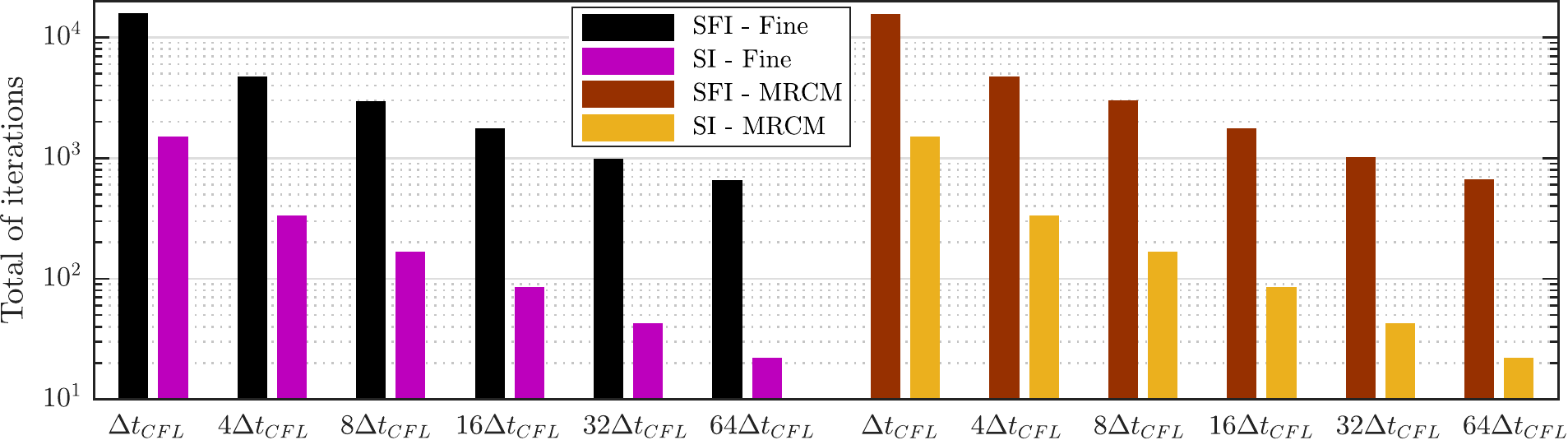}
\caption{Total accumulated of Newton iterations until time $T_{\textbf{PVI}}=0.6$ for the SI and SFI schemes combined with the fine grid velocity field and MRCM. Different sizes of time steps are considered. Note the high cost (in terms of the number of Newton iterations) of the SFI scheme when compared to the SI.}
\label{fig:Glimm_7}
\end{figure}

To close this discussion, Fig. \ref{fig:Glimm_7} shows the total accumulated of Newton iterations for approximating the 2D Riemann problem until the final time $T_{\textbf{PVI}}=0.6$. We report the total of iterations required by the SI and SFI schemes combined with the fine grid velocity field and MRCM for all the different sizes of time steps previously considered. We note the high cost (in terms of the number of Newton iterations) of the SFI scheme when compared to the SI.

Despite the fact that the SFI scheme is more expensive than the SI, the former presents more accurate solutions. We also test the SFI scheme to approximate the previous examples (in subsections \ref{A comparison of implicit solutions} and \ref{A multiscale solution}) and we only noticed very small differences between the SFI and SI approximations (in terms of saturation maps) for choices of $\Delta t\gtrsim 200 \Delta t_{CFL}$. The same accuracy for both SFI and SI schemes is observed for the mentioned examples in terms of global errors.
This is an expected result since the SI algorithm produces satisfactory approximations in the context of two-phase flows considered \cite{aziz1979petroleum}. 
However, 
for the 
problem with physics instabilities, we notice the SFI scheme performing better than the SI.  
We have shown that the MRCM works properly when combined with both SI and SFI schemes. Therefore, the MRCM combined with the SFI algorithm can be further applied to simulate more complex models.

\subsection{An example with gravity}\label{gravity}

Our last example considers an application of the sequential implicit solver to approximate a two-phase flow problem with gravity. The elliptic equation (\ref{Darcy}) when incorporating the gravity effects is given by 
\begin{equation}\label{Darcy_wg}
\begin{array}{rll}
\mathbf{u}&=-K(\mathbf{x})\big(\lambda(s)\nabla p - \lambda_g(s) g\nabla h \big) &\mbox{in}\ \Omega \\
\nabla \cdot \mathbf{u}&=q  &\mbox{in}\ \Omega \\ 
p &= p_b &\mbox{on}\ \partial\Omega_{p}\\
\mathbf{u} \cdot \mathbf{n}&= u_b, &\mbox{on}\ \partial\Omega_{u}
\end{array}
\end{equation}
where $g$ is the gravitational acceleration, $h$ is the height, and $\lambda_g$ is the gravitational mobility, that is given by
\begin{equation}
 \lambda_g(s) =\lambda_w(s)\rho_w+\lambda_o(s)\rho_o,
 \label{fluxos1_wg} 
\end{equation}
where $\rho_w$ and $\rho_o$ are, respectively, the densities of the water and oil. 
The transport equation (\ref{BL2D}) in the case with gravity is given by
\begin{equation}\label{BL2D_wg}
\begin{array}{rll}
\dfrac{\partial s}{\partial t} + \nabla \cdot \Big( f(s) \big(\mathbf{u} + K(\mathbf{x}) \lambda_o(s)(\rho_w-\rho_o)g\nabla h \big)\Big)& = 0  &\mbox{in}\ \Omega \\
s(\mathbf{x},t=0) &= s^0(\mathbf{x}) &\mbox{in}\ \Omega\\
s(\mathbf{x},t) &= \bar{s}(\mathbf{x},t). &\mbox{in}\ \partial\Omega^-
\end{array}
\end{equation}
We solve the transport problem by the Newton method as defined in Eq.  (\ref{BL2D_Newton8}). Here, the balance of fluxes $\mathcal{F}_{i,j}^{n+1,\nu}$ at a cell $(i,j)$ is a function of $f(s^{n+1,\nu})$, $\lambda_o(s^{n+1,\nu})$, $\mathbf{u}^n$, $K(\mathbf{x})$, $\rho_w$, $\rho_o$, $g$, and $\nabla h$. We use the Implicit Hybrid Upwinding (IHU) method \cite{lee2015hybrid, hamon2016implicit} to define 
\begin{equation}
\mathcal{F}_{i,j}^{n+1,\nu}=\mathcal{V}_{i,j}^{n+1,\nu}+\mathcal{G}_{i,j}^{n+1,\nu},
 \label{IHU} 
\end{equation}
where $\mathcal{V}_{i,j}^{n+1,\nu}$ represents the viscous part (i.e., the term $f(s)\mathbf{u}$ in Eq. (\ref{BL2D_wg})), while $\mathcal{G}_{i,j}^{n+1,\nu}$ represents the gravity part (i.e., the term $f(s)K(\mathbf{x}) \lambda_o(s)(\rho_w-\rho_o)g\nabla h$ 
in Eq. (\ref{BL2D_wg})). The viscous and gravity parts are treated separately in the IHU framework. 

The upwinding of the viscous term is based on the velocity field, and hence, the balance of the fluxes $\mathcal{V}_{i,j}^{n+1,\nu}$ is the same defined in Eq. (\ref{fluxes_impli}) with discrete fluxes given by Eqs. (\ref{up_fluxes_impli}) and (\ref{eq:up_fluxes_impli}).  
On the other hand, the definition of the gravity term is based on density differences, being the balance of the gravitational fluxes given by
\begin{equation}
\mathcal{G}^{n+1,\nu}_{i,j}= \tilde{G}^{n+1,\nu}_{i,j+1/2} - \tilde{G}^{n+1,\nu}_{i,j-1/2},
 \label{ihu_fluxes}
\end{equation}
where the discrete fluxes $\tilde{G}^{n+1,\nu}_{i,j\pm1/2}$ on respective interfaces $y_{j\pm1/2}$ are defined based on the fact that the heavier fluid goes down and the lighter fluid goes up as follows:
\begin{equation}
 \tilde{G}^{n+1,\nu}_{i,j\pm1/2} =
 \left\{
\begin{array}{rl}
 \Delta x\ K_{i,j\pm1/2}\dfrac{\lambda_w(s_{i,j}^{n+1,\nu})\lambda_o(s_{i,j\pm1}^{n+1,\nu})}{\lambda_w(s_{i,j}^{n+1,\nu})+\lambda_o(s_{i,j\pm1}^{n+1,\nu})}(\rho_w-\rho_o)g&\mbox{if}\ (y_{j\pm1}-y_j)g>0 \vspace{0.2cm} \\ 
  \Delta x\  K_{i,j\pm1/2}\dfrac{\lambda_w(s_{i,j\pm1}^{n+1,\nu})\lambda_o(s_{i,j}^{n+1,\nu})}{\lambda_w(s_{i,j\pm1}^{n+1,\nu})+\lambda_o(s_{i,j}^{n+1,\nu})}(\rho_w-\rho_o)g&\mbox{otherwise}
\end{array} \right. 
 \label{ihu_fluxes_impli}
\end{equation}
where $K_{i,j\pm1/2}$ is the harmonic average of $K(x_i,y_{j})$ and $K(x_i,y_{j\pm1})$. Note that the gravity effect acts along the $y$ direction.

To approximate the problem with gravity we use the extension of the trust-region algorithm with the inflection-point strategy that includes trust-regions delineated by the unit-flux and endpoints \cite{wang2013trust}. This extension also treats kinks in the Newton paths generated by the gravity term. We use an under-relaxation factor of 0.5 to ensure that the solution update does not extend the new state beyond the trust regions delineated by inflection-points and kinks.

The use of large time steps is essential for computational efficiency, especially for examples with gravitational effects that introduce restrictive CFL condition for stability, given by: 
\begin{equation}
\Delta t_{CFL}\leq \dfrac{\min \{\Delta \mathbf{x}\}}{\max|\mathfrak{F}'(s)|}, \ \text{where} \ \mathfrak{F}(s)= f(s) \big(\mathbf{u} + K(\mathbf{x}) \lambda_o(s)(\rho_w-\rho_o) g\nabla h \big).
\end{equation}
By considering the relative permeabilities $k_{ro}=(1-s)^2$ and $k_{rw}=s^2$, and hence, $f(s)=\dfrac{Ms^2}{Ms^2+(1-s)^2}$, with $M=\mu_o/\mu_w$, we have the following estimate
\begin{equation}
\begin{array}{rl}
\max|\mathfrak{F}'(s)| & \leq \max|f'(s)\mathbf{u}| + \max\big|f'(s)K(\mathbf{x}) \lambda_o(s)(\rho_w-\rho_o) g\nabla h + f(s)K(\mathbf{x}) \lambda'_o(s)(\rho_w-\rho_o) g\nabla h\big| \vspace*{0.1cm}\\
& \leq \max|f'(s)\mathbf{u}| + \left|\dfrac{(\rho_w-\rho_o)g\nabla h }{\mu_o}\right| \max\big|K(\mathbf{x})\big| \max\big| f'(s)(1-s)^2-2f(s)(1-s)\big| \vspace*{0.1cm}\\
& \leq \max|f'(s)\mathbf{u}| + 2\left|\dfrac{(\rho_w-\rho_o)g\nabla h }{\mu_o}\right| \max\big|K(\mathbf{x})\big|, \ \text{if}\ M\leq10.
\end{array}
\end{equation}
Such severe restriction highlights the importance of transport implicit methods, that allow for the use of large time steps when compared to explicit time integration approaches.

\begin{figure}
\centering
\includegraphics[scale=0.6]{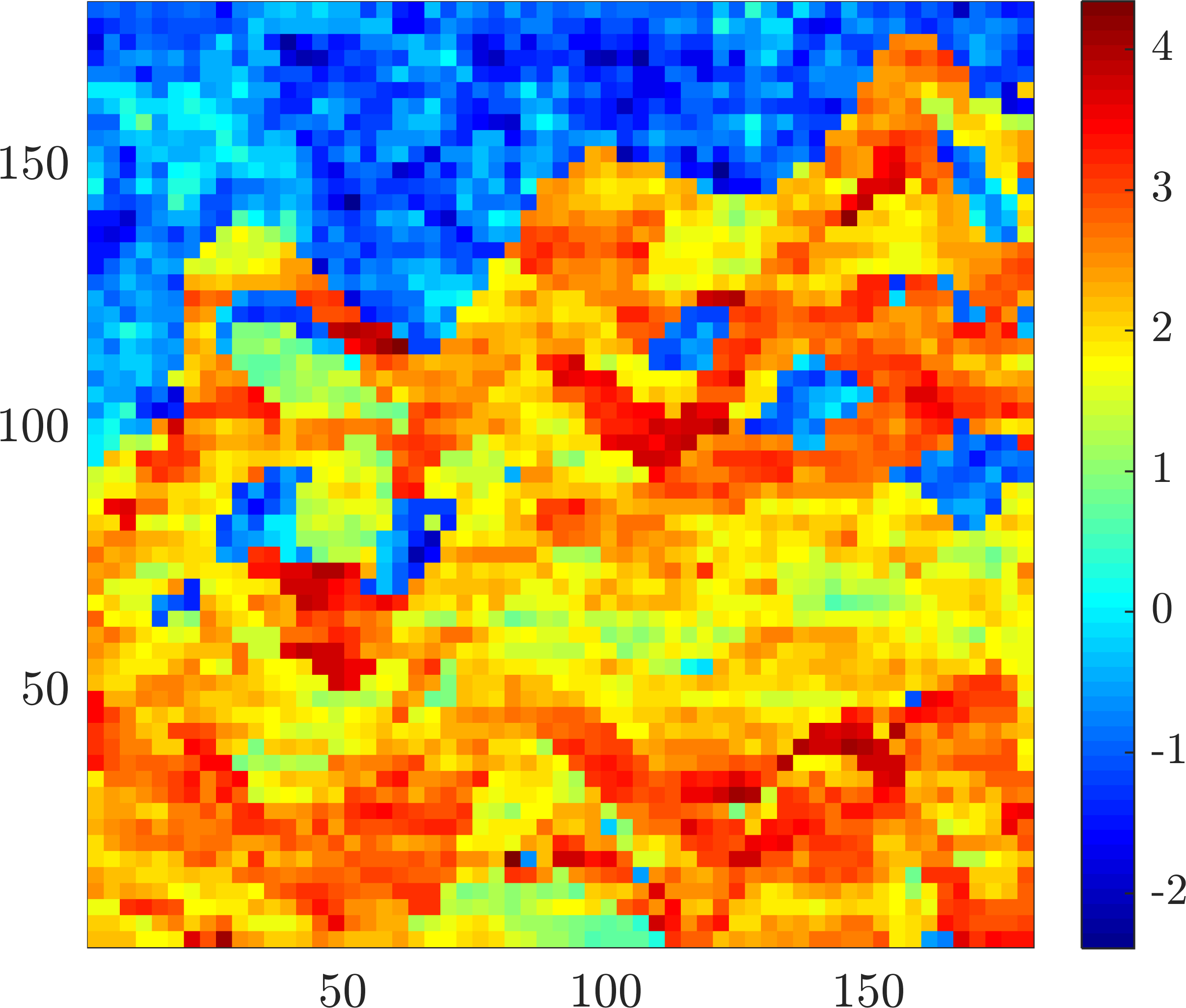}\qquad
\includegraphics[scale=0.55]{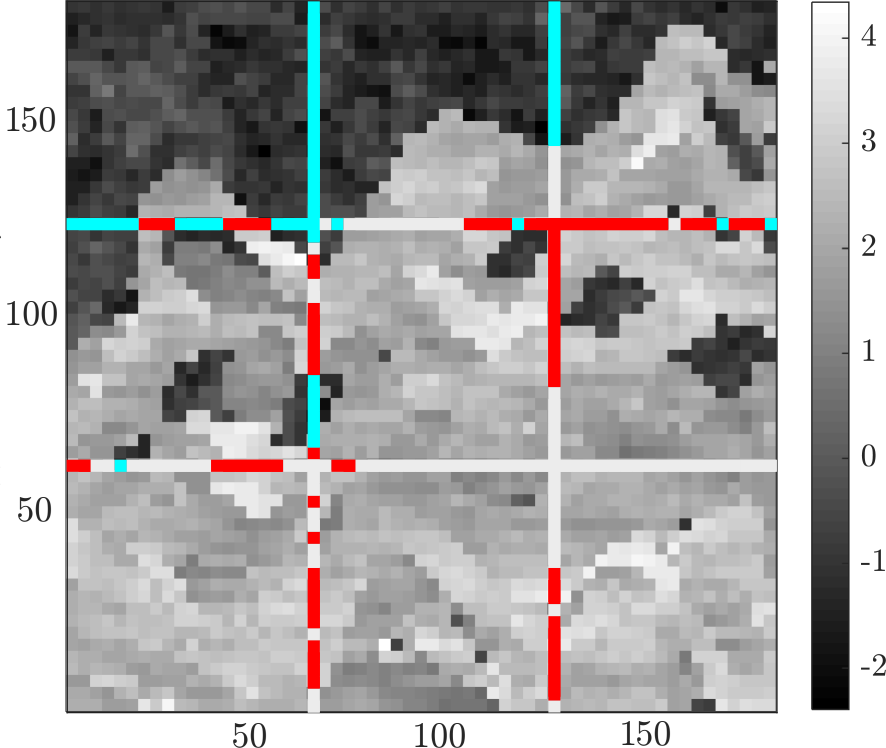}
	\caption{High-contrast permeability field (log-scaled), that is part of the top layer of the SPE10 project (left) and a map of the absolute permeability variations at the boundaries of the subdomains (right). The red and cyan colors identify the high and low-permeability regions, respectively.}
	\label{fig:perm_gravity}
\end{figure}

The geometry considered for this example is a classical quarter of a 5-spot problem \cite{chen2006computational}, where a point source for water injection is placed at the bottom left corner and a production well is placed at the top right corner of a square domain. The boundary conditions are no-flow at all boundaries. We consider the dimensionless form of the two-phase flow problem with gravity as presented in \ref{appendix}. Table \ref{table_1} shows the specification of the model, and Fig. \ref{fig:perm_gravity} (left) presents the high-contrast permeability field considered, that is part of the top layer of the SPE10 project \cite{christie2001tenth}. 
The MRCM approximation for this problem uses the adaptivity of the Robin parameter by setting $\alpha(\mathbf{x})=10^{-2}$ at the interfaces that cross highly-permeable regions and $\alpha(\mathbf{x})=10^{2}$ at the remaining interfaces. 
We denote by $a$MRCM-PBS the MRCM version that uses the adaptive Robin parameter and the physics-based interface spaces, and include for comparison the adaptive MRCM with linear interface spaces denoting by $a$MRCM-POL. In all cases a domain decomposition of $3\times 3$ subdomains with $20\times 20$ cells into each one is considered. 
Figure \ref{fig:perm_gravity} (right) shows a map of the absolute permeability variations at the boundaries of the subdomains. The red color identifies the high-permeability regions, where $\alpha(\mathbf{x})=10^{-2}$ is set and the physics-based interface spaces for pressure are defined. The cyan color identifies the low-permeability regions, where the physics-based interface spaces for flux are defined.

\begin{table}
\begin{center}
\begin{tabular}{|c c c|}
\hline 
Parameter & Value & Unit \\ 
\hline 
$L$ & 182.88 & m \\ 
Computational cells & $60\times 60$ & cells \\ 
Reference (initial) pressure & 2000 & psi \\ 
Rock permeability & Fig. \ref{fig:perm_gravity} & md \\ 
Water density & 1000 & $\mathrm{kg/m^3}$ \\ 
Oil density & 800 & $\mathrm{kg/m^3}$ \\ 
Water viscosity & 0.3 & cP \\ 
Oil viscosity & 3 & cP \\ 
Gravity acceleration & 9.80665 & $\mathrm{m/s^2}$ \\ 
Injection rate & 0.2 & PVI/y \\
\hline 
\end{tabular} 
\caption{Specification of the model for the quarter of a 5-spot problem with gravity.}
\label{table_1}
\end{center}
\end{table}

\begin{figure}
\centering
\includegraphics[scale=0.7]{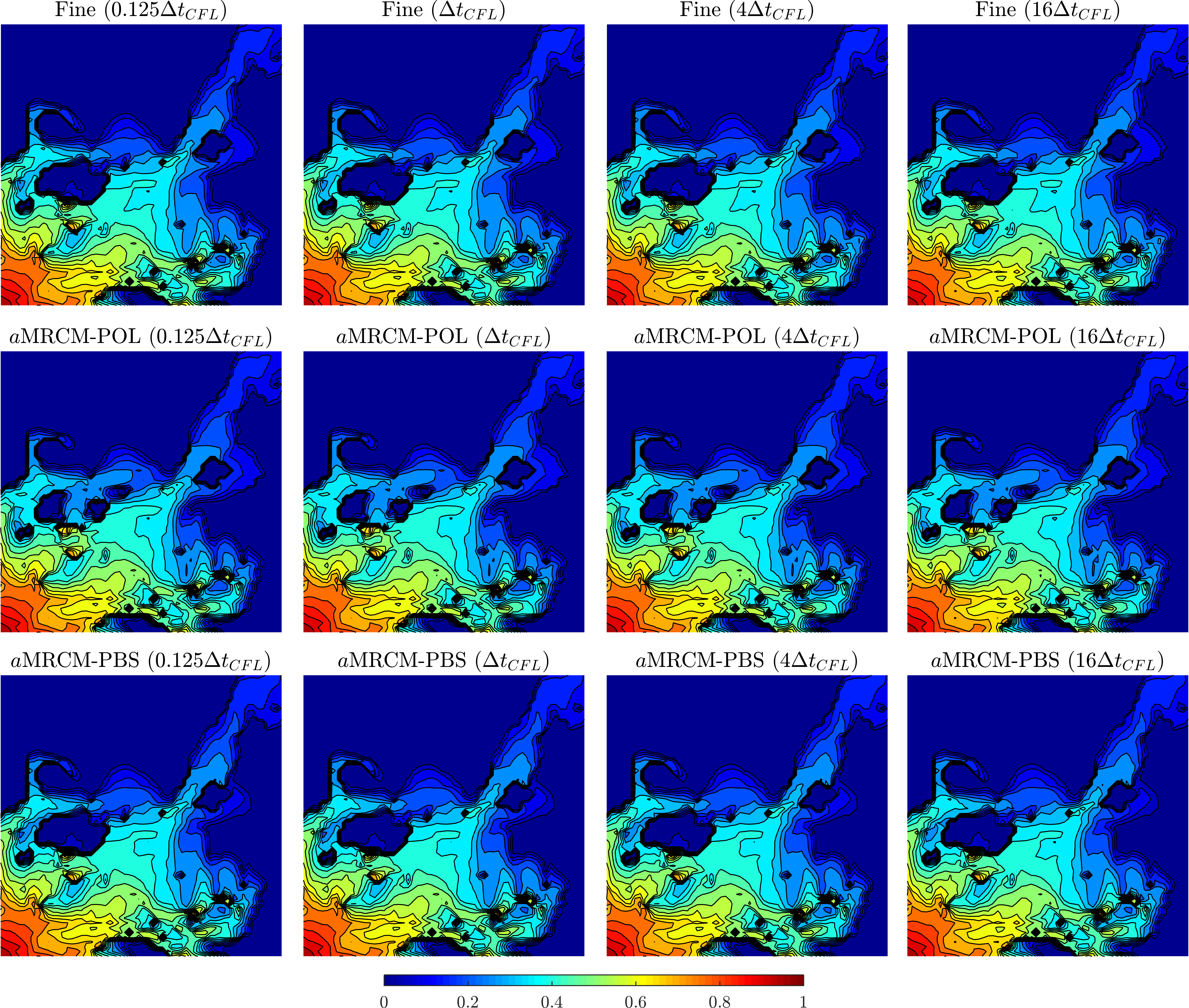}
\caption{Saturation profiles at time $T_{\textbf{PVI}}=0.21$ computed by the SI solver combined with the fine grid (top line), $a$MRCM-POL (center line), and $a$MRCM-PBS (bottom line). Different time step choices as multiples of $\Delta t_{CFL}$ are considered. All approximations are closely related.}
\label{fig:perfil_sat_gravity}
\end{figure}

Figure \ref{fig:perfil_sat_gravity} shows saturation profiles at time $T_{\textbf{PVI}}=0.21$ computed by the SI solver combined with the fine grid (top line), $a$MRCM-POL (center line), and $a$MRCM-PBS (bottom line). Different sizes of time steps are compared with the reference solutions computed with $\Delta t = 0.125\Delta t _{CFL}\approx 1.38\times10^{-5}$ (in PVI). We note that both $a$MRCM-POL and $a$MRCM-PBS produce approximations quite similar to the fine grid solution. 

\begin{figure}
\centering
\includegraphics[scale=0.65]{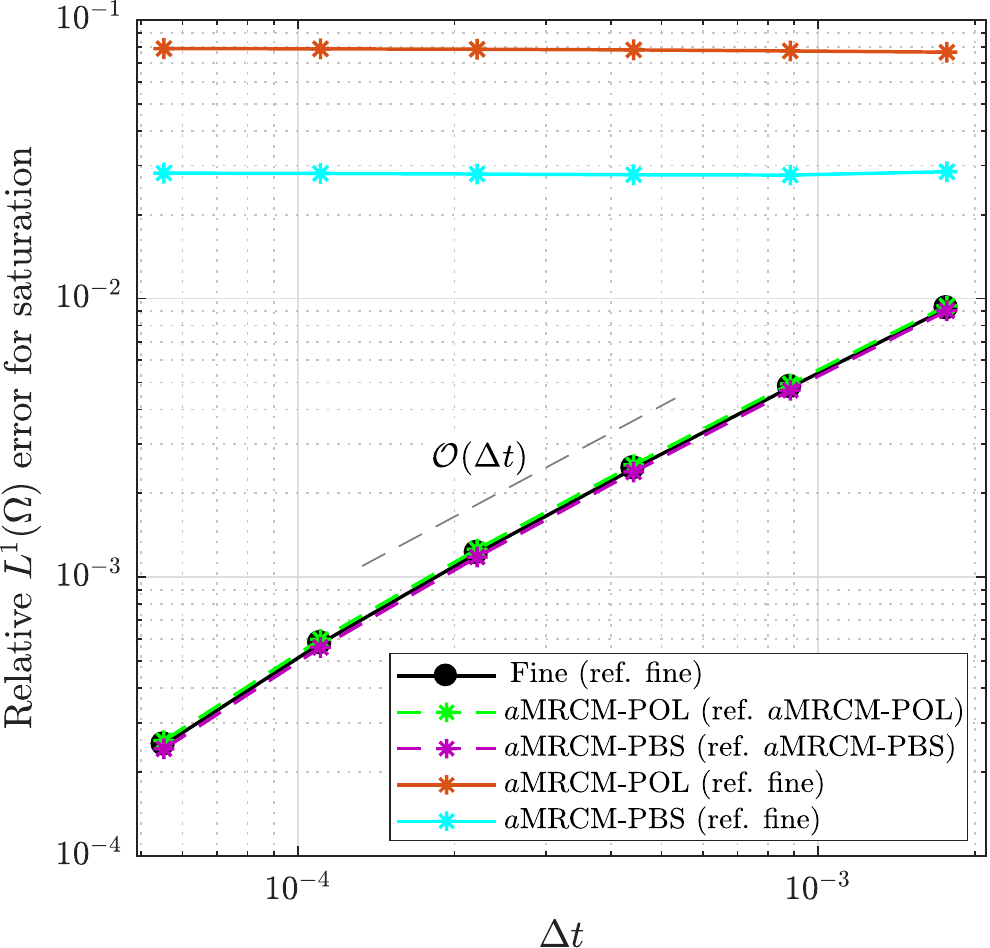}
\caption{Convergence in time of the SI scheme by setting the solution computed with $\Delta t=0.125\Delta t_{CFL}$ as reference. The fine grid, $a$MRCM-POL, and $a$MRCM-PBS are used to compute the velocity field. The errors of the $a$MRCM-POL and $a$MRCM-PBS consider their corresponding references (in space) and the fine grid reference solution. The $a$MRCM-PBS is more accurate than the $a$MRCM-POL also for problems with gravity.}
\label{fig:conv_gravity}
\end{figure}

We show in Fig. \ref{fig:conv_gravity} a convergence study for the SI scheme by setting the solution computed with $\Delta t=0.125\Delta t_{CFL}$ as reference. The fine grid, $a$MRCM-POL, and $a$MRCM-PBS are used to compute the velocity field. The errors of the $a$MRCM-POL and $a$MRCM-PBS consider their corresponding references (in space) and the fine grid reference solution. Linear slope is attained when the error of each procedure considers its corresponding spatial reference solution. The MRCM errors with respect to the fine grid solution are dominant for all time steps for both $a$MRCM-POL and $a$MRCM-PBS. The $a$MRCM-PBS improves the accuracy of the solution, reducing the errors of the $a$MRCM-POL from $8\%$ to $3\%$. Therefore the use of physics-based interface spaces is advantageous in comparison with polynomial spaces also for problems with gravity. This observation increases the importance of the $a$MRCM-PBS, which has shown accurate results, presenting error reductions up to one order of magnitude in cases with strong channelized structures \cite{rocha2020interface, rochaenhanced}. 
We remark that typical values of saturation error attained by multiscale methods are in the order of $10\%$. If this level of error from the multiscale procedure is acceptable, then very large time steps are possible making the simulation efficient. More sophisticated techniques for the multiscale methods (possibly more expensive in terms of computational cost) might be used to further reduce the errors. 

\begin{figure}
\centering
\includegraphics[scale=0.7]{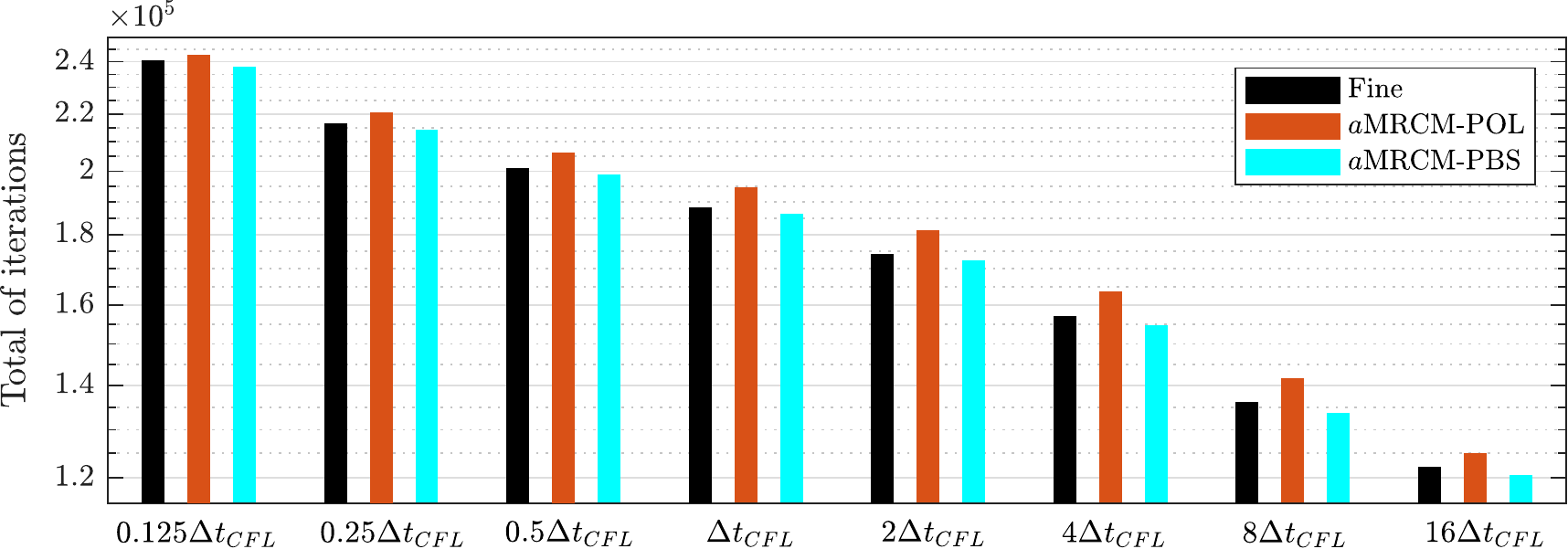}
	\caption{Total accumulated of Newton iterations until time $T_{\textbf{PVI}}=0.21$ for the SI scheme combined with the fine grid solution, $a$MRCM-POL, and $a$MRCM-PBS. Different sizes of time steps are considered. The number of iterations required by all procedures is comparable.}
	\label{fig:sum_iter_gravity}
\end{figure}

Figure \ref{fig:sum_iter_gravity} shows the total accumulated of Newton iterations until time $T_{\textbf{PVI}}=0.21$ for different sizes of time step taken as multiples of $\Delta t_{CFL}$. We start with $\Delta t=0.125\Delta t_{CFL}$ and multiply by two until $\Delta t=16\Delta t_{CFL}$. We show results for the SI scheme combined with the fine grid velocity field, $a$MRCM-POL, and $a$MRCM-PBS. The number of iterations required by all procedures is comparable (the fine mesh and $a$MRCM-PBS require essentially the same number of iterations, and the $a$MRCM-POL requires a slightly bigger number of iterations).

The MRCM produces accurate and robust results for the simulation of two-phase flows with gravity when combined with a sequential implicit scheme. Implicit schemes are fundamental to simulate practical reservoirs. 
The trust-region Newton method considered for the problem with gravity is a significant improvement to the pure Newton method \cite{wang2013trust}. However, to overcome the severe restrictions on the time step size that still appear, more specialized methods should be used, for example, the numerical trust-region solver proposed in \cite{li2015nonlinear}.

We note that the simulation converges consistently in the case of the SI scheme. On the other hand, the SFI scheme is quite inefficient for this problem, where the strong coupling between the flow and transport equations causes slow convergence of the outer-loop. One way to deal with such difficulty is presented in \cite{jiang2019nonlinear}, where the authors use nonlinear acceleration techniques to improve the convergence of the outer-loop.

\section{Conclusions} \label{sec_conclusions}

The MRCM has been combined with the SI solver for approximating two-phase flows, allowing for the use of large time steps, in contrast to conditionally-stable explicit time integration approaches. 
The numerical experiments provide strong evidence that the results produced by the MRCM combined with the SI solver are accurate and efficient. Therefore, we can replace fine grid procedures by the MRCM keeping the same parameters of the sequential implicit hyperbolic solvers. This is a promising result in terms of computational efficiency since the MRCM can take advantage of state-of-the-art parallel machines and produce two-phase flow simulations at a reduced computational cost.

To ensure convergence of the nonlinear loop, we have tested the SI solver with different trust-region algorithms. We find that the trust-region reflective and dogleg schemes are appropriate when the size of the time step is chosen of the order of $10\Delta t_{CFL}$, while the inflection-point strategy is adequate to handle sizes of time step of the order of $100\Delta t_{CFL}$ to $1000\Delta t_{CFL}$. 
We have combined the MRCM with an extension of the trust-region algorithm with the inflection-point strategy for solving problems with gravity. In the case with gravity, severe restrictions to the time step appear, making the use of transport implicit methods essential for computational efficiency. The MRCM combined with the SI scheme has shown accurate results for approximating the problem with gravitational effects. Moreover, we have also shown that the best accuracy is achieved by considering the recently introduced physics-based interface spaces for high-contrast channelized permeability fields.

The MRCM has also been combined with the SFI scheme for approximating an example with physical instabilities (fingering). We show that the MRCM works properly when combined with both SI and SFI schemes. The SFI algorithm, as well as nonlinear acceleration techniques, are currently being considered by the authors in order to further apply the MRCM to more complex flow models. 
Future works also include the implementation of the sequential implicit solver in multi-core machines.

\section*{Acknowledgments}
The authors gratefully acknowledge the financial support received from the Brazilian oil company Petrobras grant
2015/00400-4, and from the S\~ao Paulo Research Foundation FAPESP grant CEPID-CeMEAI 2013/07375-0; This study was
also financed in part by Brazilian government agencies CAPES (Finance Code 001) and CNPq grants 305599/2017-8 and 310990/2019-0.

\bibliography{bibdata}

\appendix

\section{Dimensionless form of the two-phase flow problem with gravity}\label{appendix}

We present the dimensionless form of the two-phase flow problem with gravity. Consider the two-phase flow model problem given by the following elliptic and transport equations:
\begin{equation}\label{Darcy_wg_adi}
\begin{array}{rll}
\mathbf{u}&=-K(\mathbf{x})\big(\lambda(s)\nabla p - \lambda_g(s) g\nabla h \big) &\mbox{in}\ \Omega \\
\nabla \cdot \mathbf{u}&=q  &\mbox{in}\ \Omega \\ 
p &= p_b &\mbox{on}\ \partial\Omega_{p}\\
\mathbf{u} \cdot \mathbf{n}&= u_b, &\mbox{on}\ \partial\Omega_{u}
\end{array}
\end{equation}
\begin{equation}\label{BL2D_wg_adi}
\begin{array}{rll}
\dfrac{\partial s}{\partial t} + \nabla \cdot \big( f(s) \big(\mathbf{u} + K(\mathbf{x}) \lambda_o(s)(\rho_w-\rho_o)g\nabla h \big)\big)& = 0  &\mbox{in}\ \Omega \\
s(\mathbf{x},t=0) &= s^0(\mathbf{x}) &\mbox{in}\ \Omega\\
s(\mathbf{x},t) &= \bar{s}(\mathbf{x},t) &\mbox{in}\ \partial\Omega^-
\end{array}
\end{equation}
where $ \lambda_g(s) =\lambda_w(s)\rho_w+\lambda_o(s)\rho_o$. 

In order to derive a dimensionless form of the problem (\ref{Darcy_wg_adi})-(\ref{BL2D_wg_adi}) we consider the following dimensionless quantities \newline
\begin{equation*}\label{dimensionless_1}
\mathbf{x}^*=\dfrac{\mathbf{x}}{L}, \quad 
\mathbf{u}^*=\dfrac{\mathbf{u}}{u_{\text{ref}}}, \quad
p^*=\dfrac{p}{p_{\text{ref}}}, \quad
K^*=\dfrac{K}{K_{\max}}, \quad
q^*=\dfrac{L}{u_{\text{ref}}}q, \quad
g^*=\dfrac{g}{g_{\text{ref}}}, \quad
h^*=\dfrac{h}{L},
\end{equation*}

\begin{equation*}\label{dimensionless_2}
\lambda^*=\mu_w\lambda, \qquad
\lambda^*_w=\mu_w\lambda_w, \qquad
\lambda^*_o=\mu_w\lambda_o, \qquad
\lambda^*_g=\dfrac{\mu_w}{\rho_w}\lambda_g,
\end{equation*}\newline
where $L$ is a characteristic length, $K_{\max}$ is the maximum value attained by the absolute permeability, and the reference variables $u_{\text{ref}}$, $p_{\text{ref}}$, and $g_{\text{ref}}$ are chosen such that 
\begin{equation*}
p_{\text{ref}}=\dfrac{L \mu_w u_{\text{ref}}}{K_{\max}} \quad \text{and} \quad g_{\text{ref}}=\dfrac{\mu_w u_{\text{ref}}}{\rho_w K_{\max}}.
\end{equation*}
Note that each quantity with superscript $*$ denotes a dimensionless quantity. Now, using these quantities in Eqs. (\ref{Darcy_wg_adi})-(\ref{BL2D_wg_adi}), and considering the dimensionless operator $\nabla^*=L\nabla$, we have the following dimensionless form for the elliptic equation
\begin{equation}
\begin{array}{rll}
\mathbf{u}^*&=-K^*(\mathbf{x^*})\big(\lambda^*(s)\nabla^* p^* - \lambda_g^*(s) g^*\nabla^* h^* \big) &\mbox{in}\ \Omega \\
\nabla^* \cdot \mathbf{u}^*&=q^*  &\mbox{in}\ \Omega \\ 
p^* &= p_b/p_{\text{ref}} &\mbox{on}\ \partial\Omega_{p}\\
\mathbf{u}^* \cdot \mathbf{n}&= u_b/u_{\text{ref}}, &\mbox{on}\ \partial\Omega_{u}
\end{array}
\end{equation}
and the dimensionless form for the transport equation 
\begin{equation}
\begin{array}{rll}
\dfrac{\partial s}{\partial t^*} + \nabla^* \cdot \big( f(s) \big(\mathbf{u}^* + K^*(\mathbf{x}^*) \lambda^*_o(s)\dfrac{\rho_w-\rho_o}{\rho_w}g^*\nabla^* h^* \big)\big)& = 0  &\mbox{in}\ \Omega \\
s(\mathbf{x}^*,t^*=0) &= s^0(\mathbf{x}^*) &\mbox{in}\ \Omega\\
s(\mathbf{x}^*,t^*) &= \bar{s}(\mathbf{x},t^*), &\mbox{in}\ \partial\Omega^-
\end{array}
\end{equation}
where $t^*=\dfrac{u_{\text{ref}}}{L}t$ represents the dimensionless time. However, we report the time variable in PVI (also dimensionless) in our experiments.

\end{document}